\documentclass[onefignum,onetabnum]{siamart190516}
\usepackage{theorem}
\usepackage{enumerate}
\usepackage{array}
\usepackage{amssymb}
\usepackage{mathtools}
\usepackage{latexsym}
\usepackage{makeidx}
\usepackage{fancybox}
\input epsf.sty
\usepackage{subcaption}

\newtheorem{example}[theorem]{Example}
\newtheorem{assumption}[theorem]{Assumption}
\numberwithin{equation}{section}

\outer\def\proclaim #1. #2\par{\medbreak \noindent{\bf#1.\enspace}{\sl#2}\par
  \ifdim\lastskip<\medskipamount
  \removelastskip\penalty55\medskip\fi}
\def\state #1. { \noindent{\bf#1.\enspace}}
\def\algo #1. { \noindent{\bf#1.\enspace}}

\DeclareMathOperator{\dom}{dom}

\DeclareMathOperator{\bprob}{b-prob}

\newcommand{\comp}{\,{\raise 1pt \hbox{$\scriptstyle\circ$}}\,}

\newcommand{\reals}{\mathbb{R}}
\newcommand{\Reals}{\overline{\mathbb{R}}}
\newcommand{\natnums}{{{\rm l} \kern -.13em {\rm N} }}
\newcommand{\nats}{\mathbb{N}}
\newcommand{\snats}{{I\kern -.29em N}}
\newcommand{\rats}{{Q\kern -.64em \raise 1pt \hbox{$\scriptstyle |$}\;\,}}
\newcommand{\srats}
	{{Q\kern -.56em \raise 1.2pt \hbox{$\scriptscriptstyle /$}\,}}
\newcommand{\ints}{Z\kern -.46em Z}
\newcommand{\ball}{\mathbb{B}}
\newcommand{\pluss}{\hskip1pt \raise1pt\vbox{\hrule width6pt \vskip1pt \hrule
                    width6pt} \kern-4pt{\lower1pt\hbox{\vrule height6pt
		    \kern1pt\vrule height6pt}}\hskip5pt}
\newcommand{\eop}
	{\hfill{$\vcenter{\hrule height1pt \hbox{\vrule width1pt height5pt
   	 \kern5pt \vrule width1pt} \hrule height1pt}$} \medskip}

\newcommand{\setd}{{ d \kern -.15em l}}
\newcommand{\hatsetd}{ d \hat{\kern -.15em l }}

\renewcommand{\epsilon}{\varepsilon}
\renewcommand{\phi}{\varphi}

\newcommand{\tto}{\;{\lower 1pt \hbox{$\rightarrow$}}\kern -12pt
           \hbox{\raise 2.5pt \hbox{$\rightarrow$}}\;}
\newcommand{\overto}[1]{\,{\raise 0pt\hbox{$\rightarrow$}}\kern -9pt
     \hbox{\lower 3pt \hbox{$\scriptscriptstyle#1$}}\hskip6pt}
\newcommand{\underto}[1]{\,{\lower 1pt\hbox{$\rightarrow$}}\kern -9pt
     \hbox{\raise 4pt \hbox{$\,\scriptscriptstyle#1$}}\hskip7pt}
\newcommand{\bigoverto}[1]{{\raise 0pt\hbox{$\,\longrightarrow$}}\kern -16pt
     \hbox{\lower 3pt \hbox{$\scriptscriptstyle#1$}}\hskip4pt}
\newcommand{\bigunderto}[1]{\,{\lower 1pt\hbox{$\longrightarrow$}}\kern -16pt
     \hbox{\raise 4pt \hbox{$\,\scriptscriptstyle#1$}}\hskip6pt}
\newcommand{\bigbigto}[2]{\,{\raise 0pt\hbox{$\,\longrightarrow$}}\kern -16pt
     \hbox{\lower 3pt \hbox{$\scriptscriptstyle#2$}}\kern -10pt
     \hbox{\raise 4pt \hbox{$\,\scriptscriptstyle#1$}}\hskip7pt}
\newcommand{\downto}{{\raise 1pt \hbox{$\scriptscriptstyle \,\searrow\,$}}}
\newcommand{\upto}{{\raise 1pt \hbox{$\scriptscriptstyle \,\nearrow\,$}}}

\newcommand{\notimply}
	{\quad\hbox{$\Longrightarrow \kern -14pt {/}$}\hskip6pt\quad}

\newcommand{\lto}{\,{\lower 1pt\hbox{$\rightarrow$}}\kern -10pt
     \hbox{\raise 4pt \hbox{$\, \scriptstyle l$}}\hskip7pt}
\newcommand{\eto}{\,{\lower 1pt\hbox{$\rightarrow$}}\kern -10pt
     \hbox{\raise 4pt \hbox{$\, \scriptstyle e$}}\hskip7pt}
\newcommand{\hto}{\,{\lower 1pt\hbox{$\rightarrow$}}\kern -11pt
     \hbox{\raise 4pt \hbox{$\, \scriptstyle h$}}\hskip7pt}
\newcommand{\pto}{\,{\lower 1pt\hbox{$\rightarrow$}}\kern -11pt
     \hbox{\raise 4.5pt \hbox{$\, \scriptstyle p$}}\hskip7pt}
\newcommand{\cto}{\,{\lower 1pt\hbox{$\rightarrow$}}\kern -11pt
     \hbox{\raise 4pt \hbox{$\, \scriptstyle c$}}\hskip7pt}
\newcommand{\gto}{\,{\lower 1pt\hbox{$\rightarrow$}}\kern -11pt
     \hbox{\raise 4.5pt \hbox{$\, \scriptstyle g$}}\hskip7pt}
\newcommand{\sto}{\,{\lower 1pt\hbox{$\rightarrow$}}\kern -11pt
     \hbox{\raise 4pt \hbox{$\, \scriptstyle s$}}\hskip7pt}
\newcommand{\awto}{\,{\lower 1pt\hbox{$\rightarrow$}}\kern -15pt
     \hbox{\raise 4pt \hbox{$\, \scriptstyle aw$}}\hskip7pt}
\def\Nto{\,{\raise 1pt\hbox{$\rightarrow$}}\kern -12pt
     \hbox{\lower 3pt \hbox{$\, \scriptstyle N$}}\hskip7pt}
\def\Cto{\,{\raise 1pt\hbox{$\rightarrow$}}\kern -14pt
     \hbox{\lower 3pt \hbox{$\, \scriptstyle C$}}\hskip7pt}
\def\fto{\,{\raise 1pt\hbox{$\rightarrow$}}\kern -14pt
     \hbox{\lower 3pt \hbox{$\, \scriptstyle f$}}\hskip7pt}

\newcommand{\low}[1]{{\lower1pt \hbox{$\scriptstyle #1$}}}
\newcommand{\loww}[1]{{\lower2pt \hbox{$\scriptstyle #1$}}}
\newcommand{\high}[1]{{\raise1pt \hbox{$\scriptstyle #1$}}}

\newcommand{\cB}{{\cal B}}
\newcommand{\cC}{{\cal C}}

\newcommand{\cF}{{\cal F}}

\newcommand{\cN}{{\cal N}}

\newcommand{\cR}{{\cal R}}

\newcommand{\cT}{{\cal T}}

\newcommand{\cV}{{\cal V}}

\newcommand{\nlim}{\mathop{\rm lim}\nolimits}
\newcommand{\nliminf}{\mathop{\rm liminf}\nolimits}

\newcommand{\nlimsup}{\mathop{\rm limsup}\nolimits}

\newcommand{\nnmin}{\mathop{\rm minimize}}

\newcommand{\nargmin}{\mathop{\rm argmin}\nolimits}
\newcommand{\smax}{\mathop{\rm smax}\nolimits}

\newcommand{\prob}{{\mathop{\rm prob}\nolimits}}

\newcommand{\bfxi}{\mbox{\boldmath $\xi$}}

\newcommand{\bfeta}{\mbox{\boldmath $\eta$}}

\newcommand{\bfy}{\mbox{\boldmath $y$}}

\newcommand{\lwdy}[2]{\mathrel{\mathop
        {\raisebox{0.1ex}{\null$#1$}}{\hbox{\kern -1.0em
	{\raisebox{-0.8ex}{$\scriptstyle{\;\to #2}$}}}}}}
\newcommand{\lwwdy}[2]{\mathrel{\mathop
        {\raisebox{0.2ex}{\null$#1$}}{\hbox{\kern -1.0em
	{\raisebox{-1.1ex}{$\scriptstyle{\;\to #2}$}}}}}}
\newcommand{\slwdy}[2]{\scriptsize{{\mathrel{\mathop
        {\raisebox{0.1ex}{\null$#1$}}{\hbox{\kern -1.0em
	{\raisebox{-0.8ex}{$\scriptstyle{\;\to #2}$}}}}}}}}
\newcommand{\slwwdy}[2]{\scriptsize{{\mathrel{\mathop
        {\raisebox{0.2ex}{\null$#1$}}{\hbox{\kern -1.0em
	{\raisebox{-1.1ex}{$\scriptstyle{\;\to #2}$}}}}}}}}

\definecolor{lightgray}{gray}{0.75}
\definecolor{myred}{rgb}{0.55,0,0}
\definecolor{myblue}{rgb}{0,0,0.5} % hex: #00007f
\definecolor{mygreen}{rgb}{0,0.5,0} % hex: #00007f
\definecolor{purple}{rgb}{0.5,0,0.5} % hex: #00007f
\definecolor{turq}{rgb}{0,0.805,0.816} % hex: #00007f
\definecolor{maroon}{rgb}{0.51,0,0}
\definecolor{MAROON}{rgb}{0.51,0,0}
\definecolor{redor}{rgb}{0.78,0.078,0.078}
\definecolor{dgreen}{rgb}{0,0.3,0}
\newcommand{\Ex}{\mathbb{E}}

\newcommand{\grill}{{\scriptscriptstyle\#}}
\newcommand{\bcdot}{\,{\raise .2ex \hbox{$\centerdot$}}\,}

\newcommand{\bR}{{\mathbb{R}}}
\newcommand{\bN}{{\mathbb{N}}}

\newcommand{\bP}{{\mathbb{P}}}

\newcommand{\bsxi}{{\boldsymbol{\xi}}}

\newcommand{\beq}{\begin{equation}}
\newcommand{\eeq}{\end{equation}}
\newcommand{\ba}{\begin{array}}
\newcommand{\ea}{\end{array}}

\DeclareMathOperator*{\esssup}{ess\,sup}
\DeclareMathOperator*{\essinf}{ess\,inf}
\ifpdf
\hypersetup{
	pdftitle={Performance Bounds for PDE-Constrained Optimization under Uncertainty},
	pdfauthor={P. Chen, and J.O. Royset}
}
\fi

\title{Performance Bounds for PDE-Constrained Optimization under Uncertainty
\thanks{Submitted to the editors DATE.
	\funding{The first author was partially supported by NSF under DMS-2012453, DoE under DE-SC0019303, and the Simons Foundation under award 560651. The second author was supported by the ONR Science of Autonomy
(N0001421WX00142) and AFOSR (18RT0599, 21RT0484).}}
}

\author{
	Peng Chen
	\thanks{Oden Institute for Computational Engineering and Sciences, The University of Texas at Austin, Austin, TX 78712. Email: \email{peng@oden.utexas.edu}.}
\and Johannes O. Royset
\thanks{Operations Research Department, Naval Postgraduate School, Monterey, CA 93943. Email: \email{joroyset@nps.edu}.}
}

\begin{document}
	
	\maketitle
	
	% REQUIRED
	\begin{abstract}
Computational approaches to PDE-constrained optimization under uncertainty may involve finite-dimensional approximations of control and state spaces, sample average approximations of measures of risk and reliability, smooth approximations of nonsmooth functions, penalty approximations of constraints as well as many other kinds of inaccuracies. In this paper, we analyze the performance of controls obtained by an approximation-based algorithm and in the process develop estimates of optimality gaps for general optimization problems defined on metric spaces. Under mild assumptions, we establish that limiting controls have arbitrarily small optimality gaps provided that the inaccuracies in the various approximations vanish. We carry out the analysis for a broad class of problems with multiple expectation, risk, and reliability functions involving PDE solutions and appearing in objective as well as constraint expressions. In particular, we address problems with buffered failure probability constraints approximated via an augmented Lagrangian.  We demonstrate the framework on an elliptic PDE with a random coefficient field and a distributed control function.
	\end{abstract}
	
	% REQUIRED

	% REQUIRED

\section{Introduction}

Optimization problems involving partial differential equations (PDEs) arise widely in design and control of physical systems \cite{HinzePinnauUlbrichEtAl09, Troeltzsch10}. Recent applications include shape optimization \cite{KolvenbachLassUlbrich18, GeiersbachLoayza-RomeroWelker20}, turbulent combustion \cite{ChenVillaGhattas19}, acoustic wave propagation \cite{YangGunzburger17}, metamaterials \cite{ChenHabermanGhattas21}, and plasma fusion \cite{WechsungGiulianiLandremanEtAl21}. Problems of this kind usually come with significant uncertainty in the form of unknown external loadings, material coefficients, boundary or initial conditions, geometries, and other factors, which leads to PDE-constrained optimization problems under uncertainty. Formulations of these problems leverage scalarizations of the random quantities of interest via expectations, variances \cite{AlexanderianPetraStadlerEtAl17}, superquantiles (a.k.a. conditional/average value-at-risk) \cite{KouriSurowiec16, GarreisSurowiecUlbrich21}, worst-case measures of risk \cite{LassUlbrich17}, and probabilities of failure \cite{ChenGhattas21}.
These formulations tend to involve integration with respect to a probability measure, which needs to be approximated by numerical integration, such as Monte Carlo and quasi Monte Carlo techniques \cite{RoemischSurowiec21,MartinKrumscheidNobile21,GuthKaarniojaKuoEtAl19}, sparse grid quadrature methods \cite{ChenQuarteroniRozza13}, or variance-reduction techniques based on Taylor expansions \cite{ChenVillaGhattas19}. These approximations compound the already significant challenges associated with discretization and solution of the underlying PDE. Techniques to mitigate the overall computational cost include those based on model reduction \cite{ChenQuarteroni14,LassUlbrich17}, multilevel \cite{ChenQuarteroniRozza16,AliUllmannHinze17} and multifidelity \cite{NgWillcox14} approximations,  low-rank tensor decomposition \cite{BennerOnwuntaStoll16}, stochastic collocation \cite{KouriHeinkenschlossRidzalEtAl13}, stochastic Galerkin \cite{KunothSchwab16}, and Taylor approximation \cite{AlexanderianPetraStadlerEtAl17}.
In this paper, we analyze a multitude of approximations in PDE-constrained optimization under uncertainty. We focus on a broad class of problems involving several quantities of interest and, thus, several expectation, risk, and reliability functions. These quantities may arise in an objective function, as constraint functions, or both.

In addition to discretization of state and control spaces and numerical evaluation of integrals, the broad class of problems gives rise to several other approximations. Smooth functions approximate nonsmooth ones to facilitate the use of gradient-type algorithms. Penalty and Lagrangian terms model constraints. Weights aggregate multiple objective functions. The combined effect of this multitude of approximations on the obtained solutions is nontrivial and increasingly complex as we move towards more sophisticated optimization problems involving novel performance criteria and numerous quantities of interest. We lay out conditions under which small errors in the various approximations indeed lead to a guaranteed good performance of a control computed using these approximations. This does not hold in general. It is well known that approximating functions may converge pointwise to an actual function but have minimizers and minima far from those of the actual problem \cite[Figure 4.5]{primer}. More concretely, a convergent Runga-Kutta method employed to solve an ordinary differential equation might still induce errors that prevent minimizers of the approximating optimization problems from converging to a minimizer of the actual problem \cite{SchwartzPolak.96,Hager.00}. We present a unified framework that addresses a wide array of approximations and applications. It stretches beyond PDE-constrained optimization, but also specializes to an example with a buffered failure probability constraint and to problems defined in terms of an elliptic PDE with a random coefficient field and a distributed control function.

The various inaccuracies lead to approximating optimization problems that subsequently need to be solved by optimization algorithms. Since these problems are rarely convex, one cannot expect their minimizers to be within reach. We furnish performance bounds and optimality gaps for {\em any} solution---local, global, or neither---obtained from the approximating problems. Still, better solutions of the approximating problems naturally translate into tighter bounds.

The reasoning towards performance bounds and optimality gaps follows a novel breakdown of errors into those caused by discretization of the control space on one side and those produced by all other approximations on the other side. This division is motivated by the fact that the former imposes a restriction on the problem while the latter can swing either way. As a result, we often analyze the ``other approximations'' in a setting of finite dimensions and this reduces the required assumptions. We omit an analysis of optimality conditions and the assumptions required to ensure convergence of stationary points of approximating problems to those of the actual problem; see for example \cite{PhelpsRoysetGong16, KouriSurowiec18} for efforts in this direction.
Our results provide the foundation for numerous algorithms, including those involving adaptive refinements. However, we defer detailed algorithms as well as numerical demonstrations to future publications.

Our main technical tool is epi-convergence; see for example \cite[Chapter 7]{vaan}, \cite[Section 4.C]{primer}, and references therein. It is well known that epi-convergence of approximating functions ensures that the corresponding minimizers can only converge to a minimizer of the limiting function. However, epi-convergence can be difficulty to verify or might simply not hold for realistic problems involving a multitude of approximations. Many studies concentrate on a single source of approximation;  \cite{KouriSurowiec20} deals with nonsmoothness and \cite{Diem,SurBirge} with sampling. The challenge of combining discretization of an underlying space with sample averages is illustrated by the less-than-ideal convergence results of Theorem 3.14 in \cite{RoysetWets.20}, which anyhow fails to address other approximations and expectation constraints. The difficulties compound for PDE-constrained optimization problems because the approximating problems solved by numerical algorithms are defined on finite-dimensional Euclidean spaces while the actual problem resides in infinite dimensions.

The rest of the paper is organized as follows: Section \ref{sec:formulation} considers formulations and illustrative examples. Section \ref{sec:general-approximation} develops a general approximation theory and, thus, lays the foundation for addressing PDE-constrained optimization problems in Section \ref{sec:performance-bounds}. Section \ref{sec:applications} demonstrates the theory in the context of optimal control of an elliptic PDE with a random coefficient field and two quantities of interest.

\section{Problem Formulations}
\label{sec:formulation}

We consider the stochastic optimization problem
\begin{equation}\label{eqn:actualpde0}
	\nnmin_{z\in Z} ~\phi(z) = \iota_{A}(z) + f(z) + h\Big( \Ex\big[G(\bfxi,z)\big] \Big),
\end{equation}
referred to as the {\em actual problem}, where the {\em control} $z$ resides in a separable Banach space $(Z,\|\cdot\|_Z)$ and is restricted to an {\em admissible set} $A\subset Z$ as expressed by the {\em indicator function}; for any set $C$, $\iota_C(c) = 0$ if $c \in C$ and $\iota_C(c) = \infty$ otherwise. The {\em objective function} $\phi:Z\to \Reals = [-\infty, \infty]$ is further given by a {\em cost function} $f:Z\to \reals$. Main complications stem from the last term in the objective function, which is defined by a {\em random field} $\bfxi$, an $m$-dimensional vector of {\em quantities of interest} $G(\bfxi,z)$, and a {\em monitoring function} $h:\reals^{m}\to \Reals$. The monitoring function assesses the various quantities of interest, which in turn depend implicity on the solution of a {\em PDE} defined on a domain $D\subset\reals^d$ and parametrized by $z\in Z$ and realizations of $\bfxi$.
	
Formally, let $(\Omega,\cF, \mathbb{P})$ be a probability space and $\Xi$ be a Hilbert space of functions on $D$ equipped with the Borel $\sigma$-field $\cB$. On the probability space, we define a random field, i.e., a measurable function $\bfxi:\Omega\to \Xi$. We denote by $\xi$ a realization of $\bfxi$, i.e., $\xi = \bfxi(\omega) \in \Xi$ for some $\omega \in \Omega$. Thus, boldface letters indicate a random element and regular font its realization. In the usual way, $\bfxi$ defines another probability space $(\Xi,\cB,P)$, where $P(B) = \mathbb{P}(\{\omega\in \Omega~|~\bfxi(\omega) \in B \})$ for $B\in \cB$. Without imposing any practical limitation, we let this probability space be complete.
	
We assume that the solutions of the PDE lie in a separable Banach space $(U,\|\cdot\|_U)$ and, for given $\xi\in \Xi$ and $z\in Z$, the PDE has a unique solution in $U$ denoted by $s(\xi,z)$. Section \ref{sec:applications} furnishes details about the existence of such solutions. For the majority of the development, it suffices to recognize that the solutions of the PDE for various $\xi$ and $z$ are given by a {\em solution mapping} $s:\Xi\times Z\to U$.
	
The quantities of interest are defined by the {\em performance functions} $g_i:U\times Z\to \reals$, $i=1, \dots, m$, and the solution mapping $s$:
\begin{equation}\label{eq:performance-fun}
G:\Xi\times Z\to \reals^m, ~~\mbox{ with } G(\xi,z) = \Big(g_1\big(s(\xi,z),z\big), \dots, g_m\big(s(\xi,z),z\big)\Big).
\end{equation}
For any $z\in Z$, we write $\Ex\big[G(\bfxi,z)\big] = \int G(\xi,z) dP(\xi)$, which is assumed to be well-defined and finite as discussed below.
	
We adopt the usual rules for extended arithmetic. In particular, $\infty - \infty = \infty$; see \cite[Section 1.D]{primer}. Thus, $z\not\in A$ and/or $h(\Ex[G(\bfxi,z)])=\infty$ produce $\phi(z) = \infty$, which implies that $z$ is infeasible. With the convention $\nargmin \phi = \{z\in Z~|~\phi(z) = \inf \phi < \infty\}$, where $\inf \phi = \inf\{\phi(z)~|~z\in Z\}$, this set of minimizers cannot include such $z$.

\subsection{Illustrative Examples}
\label{sec:examples}

The actual problem \eqref{eqn:actualpde0} addresses a vast array of applications. We highlight some possibilities.
	
\begin{example}{\rm (expectation minimization with regularization).}\label{eExpectationReg} Suppose that we seek to determine a control $z\in A \subset Z = L^2(D;\reals^q)$, the space of Lebesgue square-integrable functions from $D$ to $\reals^q$, that minimizes the expected system performance as quantified by $g_1:U\times Z \to \reals$. This leads to the problem
\[
\nnmin_{z\in A} ~\Ex\Big[g_1\big(s(\bfxi,z),z\big)\Big] + \theta\|z\|^2_Z,
\]
where a regularization term is included with $\theta \in [0,\infty)$. The problem is of the form \eqref{eqn:actualpde0}: set $m=1$, $h(w) = w$ and $f(z) = \theta\|z\|_Z^2$.
\end{example}

\begin{example}{\rm (tracking objective under expectation constraint).} While using the same admissible controls as in Example \ref{eExpectationReg}, suppose that we seek to bring the random state $s(\bfxi,z)$ as close as possible to a target state $\bar u\in U$ on average and force a performance function $g_1:U\times L^2(D;\reals^q) \to \reals$ to lie between $\alpha$ and $\beta$ on average. This leads to the problem
\[
\nnmin_{z\in A} ~\Ex\Big[ \big\| s(\bfxi,z) - \bar u\big\|_{U}\Big]  ~\mbox{ subject to } ~ \alpha \leq \Ex\Big[g_1\big(s(\bfxi,z),z\big)\Big] \leq \beta,				
\]
which is of the form \eqref{eqn:actualpde0} with $m=2$, $h(w) = \iota_{[\alpha,\beta]}(w_1)+ w_2$, $f(z) = 0$, and $g_2(u,z) = \| u - \bar u\|_{U}$. The bounds $\alpha\leq \beta$ could coincide to produce an equality constraint or they might have $\alpha = -\infty$ and $\beta = 0$ to produce the inequality $\Ex[g_1(s(\bfxi,z),z)] \leq 0$.
\end{example}

\begin{example}{\rm (risk modeling).}\label{eRiskMin} Modeling focused on ``worst-case'' performance instead of ``average'' performance is accomplished using a regular measure of risk  that maps random variables with finite second moments into $\Reals$ \cite[Chapter 8]{RockafellarRoyset.15,primer}. Let $\cR_i$, $i=1, \dots, m$, be a collection of regular measures of risk. Suppose that we seek a control $\hat z\in \hat A\subset \hat Z$ that minimizes the risk associated with a performance function $\hat g_1:U\times \hat Z\to \reals$ and satisfies constraints of nonnegative risk related to $\hat g_i:U\times \hat Z\to \reals$, $i=2, 3, \dots, m$. This leads to the problem
\[
\nnmin_{\hat z\in \hat A} ~\cR_1\Big( \hat g_1\big(s(\bfxi,\hat z),\hat z\big) \Big) ~\mbox{ subject to } ~ \cR_i \Big(\hat g_i\big(s(\bfxi,\hat z),\hat z\big)\Big) \leq 0,~~i=2, 3, \dots, m,
\]
which turns out to be expressible in the form \eqref{eqn:actualpde0} for common measures of risk.
\end{example}
\state Detail. By \cite[Thm. 8.9]{primer}, every regular measure of risk $\cR$ can be expressed as $\cR(\bfeta) = \min_{\gamma\in\reals} \gamma + \cV(\bfeta-\gamma)$ for some regular measure of regret $\cV$, which also maps random variables with finite second moments into $\Reals$. Thus, the problem is equivalently stated as a minimization problem over $\hat z$ and auxiliary variables $\gamma_1, \dots, \gamma_m$ using regular measures of regret $\cV_1, \dots, \cV_m$ corresponding to $\cR_1, \dots, \cR_m$, respectively. Many common measures of regret $\cV$ are of the expectation kind, which means that for some
convex function $v:\reals\to \reals$ one has $\cV(\bfeta) = \Ex[v(\bfeta)]$ for random variables $\bfeta$; see \cite{RockafellarUryasev.13} and \cite[Chapter 8]{primer}. For $\alpha \in (0,1)$, the {\em penalty regret} $\cV(\bfeta) = \Ex[\max\{0,\bfeta\}]/(1-\alpha)$ is a prominent example producing the superquantile risk (a.k.a. CVaR and AVaR) \cite{RockafellarUryasev.00}. If $\cV_1, \dots, \cV_m$ are of the expectation kind and expressible using the convex functions $v_1, \dots, v_m$, then the problem takes the equivalent form
\begin{align*}
& \nnmin_{\hat z\in \hat A, \gamma_1, \dots, \gamma_m} ~\Ex\Big[\gamma_1 + v_1\Big( \hat g_1\big(s(\bfxi,\hat z),\hat z\big) -\gamma_1\Big)\Big]\\
& ~~\mbox{ subject to } ~~ \Ex\Big[\gamma_i + v_i \Big(\hat g_i\big(s(\bfxi,\hat z),\hat z\big)-\gamma_i\Big)\Big] \leq 0,~~~i=2, 3, \dots, m,
\end{align*}
which indeed is of the form \eqref{eqn:actualpde0} with control space $Z = \hat Z \times \reals^m$, set of admissible controls $A = \hat A \times \reals^m$, performance functions
\[
g_i(u,z) = \gamma_i + v_i\big(\hat g_i(u,\hat z)-\gamma_i\big), ~i=1, \dots, m, ~~ \mbox{ for } ~~z = (\hat z, \gamma_1, \dots, \gamma_m)
\]
and monitoring function $h(w) = w_1 + \iota_{(-\infty, 0]}(w_2) + \cdots + \iota_{(-\infty, 0]}(w_m)$.\eop

\begin{example}{\rm (buffered failure probability constraint).}\label{eBuffered} A reliability constraint can beneficially be expressed in terms of a buffered failure probability \cite{RockafellarRoyset.10}, which is better behaved mathematically than a failure probability \cite{MafusalovUryasev.18,ByunRoyset.21}. For a random variable $\bfeta$ with finite mean, the buffered failure probability is defined as
\[
\bprob\{\bfeta > 0\} = \begin{cases}
0 & \mbox{ if } ~\prob\{ \bfeta > 0\} = 0\\
1-\alpha & \mbox{ if } ~\prob\{ \bfeta > 0\} > 0~\mbox{ and } ~\Ex[\bfeta] < 0\\
1 & \mbox{ otherwise,}
\end{cases}
\]
where $\alpha\in (0,1)$ is the probability that makes the $\alpha$-superquantile $\bar Q_\alpha(\bfeta)$ of $\bfeta$ equal to zero. As expressed by the $\alpha$-quantile $Q_\alpha(\bfeta)$ of $\bfeta$, the $\alpha$-superquantile is defined as \cite[Section 3.C]{primer}:
\[
\bar Q_\alpha(\bfeta) = Q_\alpha(\bfeta) + \frac{1}{1-\alpha} \Ex\Big[\max\big\{ 0, \bfeta - Q_\alpha(\bfeta)\big\}\Big].
\]
The problem of determining a control $\hat z\in \hat A \subset \hat Z$ that minimizes a cost function $\hat f:\hat Z\to \reals$ plus the average value of a quantity of interest represented by $\hat g_1:U\times \hat Z \to \reals$ and produces sufficient reliability relative to another quantity of interest expressed by $\hat g_2:U\times \hat Z \to \reals$ is formulated as
\begin{equation}\label{eqn:bufferProbForm0}
\begin{split}
&\nnmin_{\hat z\in \hat A} \hat f(\hat z) + \Ex\Big[\hat g_1\big(s(\bfxi,\hat z),\hat z\big)\Big]\\
&\mbox{ subject to } \bprob\Big\{\hat g_2\big(s(\bfxi,\hat z),\hat z\big)>0\Big\} \leq 1-\alpha.
\end{split}
\end{equation}
The problem can be stated in the form \eqref{eqn:actualpde0}.
\end{example}
\state Detail. For $\alpha \in (0,1)$, $\bar Q_\alpha(\bfeta) \leq 0$ if and only if $\bprob\{ \bfeta > 0\} \leq 1-\alpha$ \cite[Section 3.E]{primer}. This fact together with the discussion in Example \ref{eRiskMin} about the penalty regret imply that \eqref{eqn:bufferProbForm0} is equivalent to
\begin{equation}\label{eqn:bufferProbForm}
\begin{split}
& \nnmin_{\hat z\in \hat A, \gamma\in\reals} \hat f(\hat z) + \Ex\Big[\hat g_1\big(s(\bfxi,\hat z),\hat z\big)\Big] \\
& \mbox{ subject to } ~~ \Ex\bigg[ \gamma + \frac{1}{1-\alpha} \max\Big\{0, \hat g_2\big(s(\bfxi,\hat z),\hat z\big)-\gamma\Big\} \bigg] \leq 0.
\end{split}
\end{equation}
This problem is of the form \eqref{eqn:actualpde0} with $Z = \hat Z \times \reals$, $A = \hat A \times \reals$,  $h(w_1,w_2) = w_1 + \iota_{(-\infty, 0]}(w_2)$ and, for  $z = (\hat z, \gamma)$, $f(z) = \hat f(\hat z)$, $g_1(u, z) = \hat g_1(u,\hat z)$, and $g_2(u, z) = \gamma + \max\{0, \hat g_2(u,\hat z)-\gamma\}/(1-\alpha)$.\eop

\subsection{Approximations}\label{subsec:approx}
	
The examples illustrate the breadth of situations addressed by \eqref{eqn:actualpde0}, but also the significant challenges associated with the development and justification of computational procedures. In addition to the need for discretization of the control space $Z$ and numerical solution of the underlying PDE, $m$ expectations must be estimated and potentially nonsmooth performance functions might emerge from reformulations of risk measures (see Examples \ref{eRiskMin} and \ref{eBuffered}) causing the need for smoothing.  Nonlinear constraints commonly require approximations via augmented Lagrangian and penalty formulations. These considerations produce approximations of \eqref{eqn:actualpde0}, where $Z$ is replaced by a finite-dimensional subspace, the monitoring function $h$ is replaced by smooth alternatives, expectations are replaced by sample average approximations and the quantities of interest, as summarized by $G$, are replaced by approximations capturing numerical solutions of the underlying PDE, and also smooth approximations necessitated by nonsmooth performance functions.

In the setting of Example \ref{eBuffered}, one might face five types of approximations: A finite-dimensional subspace $\hat Z^n$ restricts $\hat Z$, a sample $\{\xi_1, \dots, \xi_\nu\}$ furnishes an estimate of expectations, an approximating solution mapping $s^\nu$ replaces $s$, an augmented Lagrangian term with parameter $\theta^\nu \in [0,\infty)$ and multiplier $y^\nu\in\reals$ substitute for the constraint, and $\smax(\gamma;\beta^\nu) = \beta^\nu \ln (1+\exp(\gamma/\beta^\nu))$ approximates the max-function $\gamma\mapsto \max\{0,\gamma\}$ with an error of at most $2\beta^\nu$ using a tunable parameter $\beta^\nu \in (0,\infty)$; see \cite[Ex. 4.16]{primer}.

Before implementing these approximations, we reformulate \eqref{eqn:bufferProbForm} using a nonnegative slack variable $\sigma$:
\begin{equation}\label{eqn:bufferProbFormLagr}
\begin{split}
& \nnmin_{\hat z\in \hat A, \gamma\in\reals, \sigma \geq 0} \hat f(\hat z) + \Ex\Big[\hat g_1\big(s(\bfxi,\hat z),\hat z\big)\Big]\\
& \mbox{ subject to } ~~ \Ex\bigg[ \sigma + \gamma + \frac{1}{1-\alpha} \max\Big\{0, \hat g_2\big(s(\bfxi,\hat z),\hat z\big)-\gamma\Big\} \bigg]  = 0.
\end{split}
\end{equation}
The problem is of the form \eqref{eqn:actualpde0} with $Z = \hat Z \times \reals^2$, $A = \hat A \times \reals \times [0,\infty)$, $h(w_1,w_2) = w_1 + \iota_{\{0\}}(w_2)$ and, for  $z = (\hat z, \gamma, \sigma)$, one has $f(z) = \hat f(\hat z)$, $g_1(u, z) = \hat g_1(u,\hat z)$ and
\begin{equation}\label{eq:g2uz}
	g_2(u, z) = \sigma + \gamma + \frac{1}{1-\alpha}\max\big\{0, \hat g_2(u,\hat z)-\gamma\big\}.	
\end{equation}
An approximating problem for \eqref{eqn:bufferProbFormLagr} is then
\begin{equation}\label{eqn:bufferProbFormLagrApprox}
\begin{split}
		& \nnmin_{\hat z\in \hat A\cap \hat Z^n, \gamma\in\reals, \sigma\geq 0} \hat f(\hat z) + \frac{1}{\nu} \sum_{j=1}^\nu \hat g_1\big(s^\nu(\xi_j,\hat z),\hat z\big)\\
		& ~~~~~~~~~+ \frac{y^\nu}{\nu} \sum_{j=1}^\nu \sigma + \gamma + \frac{1}{1-\alpha} \smax\Big(\hat g_2\big(s^\nu(\xi_j,\hat z),\hat z\big)-\gamma; ~\beta^\nu\Big)\\
		& ~~~~~~~~~+ \theta^\nu \Bigg(\frac{1}{\nu} \sum_{j=1}^\nu \sigma + \gamma + \frac{1}{1-\alpha} \smax \Big(\hat g_2\big(s^\nu(\xi_j,\hat z),\hat z\big)-\gamma;~\beta^\nu\Big)\Bigg)^2,
\end{split}
\end{equation}
where the objective function being minimized is an augmented Lagrangian \cite[Ex. 6.10]{primer}. Presumably, the approximating problem does not involve any difficult constraints, has a smooth objective function in many practical settings and can be addressed using standard nonlinear programming solvers.

One could hope that solutions of the approximating problem produce reasonable solutions for \eqref{eqn:bufferProbForm0} when the subspace $\hat Z^n$ is sufficiently close to $\hat Z$, the sample size $\nu$ is large, the numerical solution mapping $s^\nu$ approximates $s$ rather well, the penalty parameter $\theta^\nu$ is high, the multiplier $y^\nu$ is suitably selected, and the smoothing parameter $\beta^\nu$ is near zero. Of course, this cannot be expected in general. The next sections examine performance bounds for such approximations in the concrete setting of \eqref{eqn:bufferProbFormLagr} and also more broadly.

\section{General Approximation Theory}	\label{sec:general-approximation}
	
To prepare a foundation for a unified treatment of the actual problem \eqref{eqn:actualpde0}, we approach an abstract optimization problem and its approximations. Suppose that $(X,d_X)$ is a metric space\footnote{The metric space setting here hints to the possibility of extending the reach of later results beyond Banach spaces.} and let $\nats = \{1, 2, \dots\}$. For $f:X\to \Reals$, consider the problem
	\begin{equation}\label{eqn:actualpr}
		\nnmin_{x\in X} f(x)
	\end{equation}
	and the discretized problems
	\begin{equation}\label{eqn:discr}
		\Big\{\nnmin_{x\in X} f(x) + \iota_{X^n}(x), ~~n\in \nats\Big\},
	\end{equation}
	where $X^n\subset X$ corresponds to a finite-dimensional space. Specifically, for each $n\in\nats$, we assume there is a surjective mapping $T_n:\reals^n\to  X^n$, i.e., for every $x\in X^n$, one can identify a vector $x_n \in \reals^n$ satisfying $x = T_n(x_n)$.  Through this correspondence, we identify the finite-dimensional problems
	\begin{equation}\label{eqn:discrF}
		\Big\{\nnmin_{x_n\in \reals^{n}} f_n(x_n), ~~n\in\nats\Big\},
	\end{equation}
	where $f_n:\reals^{n}\to \Reals$ is given by
	\begin{equation}\label{eqn:fn}
		f_n(x_n) = f\big( T_n(x_n)\big).
	\end{equation}
	We see that \eqref{eqn:discr} and \eqref{eqn:discrF} are equivalent, with
	\begin{equation}\label{eqn:equivInf}
		\inf\big\{f_n(x_n)~\big|~ x_n\in \reals^{n}\big\} =  \inf\big\{f(x) + \iota_{X^n}(x)~\big|~ x\in X\big\}.
	\end{equation}
	
	Although finite dimensional, \eqref{eqn:discrF} might require approximations for computational and other reasons. We consider the collection of approximating problems
	\begin{equation}\label{eqn:approx}
		\forall n\in\nats: ~~~\Big\{\nnmin_{x_n\in \reals^{n}} f_n^\nu(x_n), ~~\nu\in\nats\Big\},
	\end{equation}
	where $f_n^\nu:\reals^{n}\to \Reals$ is an approximation of $f_n$ that presumably becomes more accurate as $\nu\to \infty$. We seek to justify the following approach to solving \eqref{eqn:actualpr}: apply an optimization algorithm to one or more of the approximating problems \eqref{eqn:approx} and obtain a solution that serves as an approximating solution of \eqref{eqn:actualpr}. As a preview of our treatment of the actual problem \eqref{eqn:actualpde0}, the finite-dimensional problems \eqref{eqn:discrF} correspond to a discretization of the control space and the approximating problems \eqref{eqn:approx} implement all other approximations as well, including discretization adopted to solve the underlying PDE. Since discretization of the control space can be viewed as just another approximation, it appears somewhat arbitrary to single out that discretization. However, this separation is technically convenient and allows us to utilize the fact that discretization of the control space imposes a restriction on the problem while the other approximations may not have that characteristic.
	
	We adopt the following notation and concepts. For $\delta \in [0,\infty)$ and a function $g:X\to \Reals$, the set of near-minimizers is given by
	\[
	\delta \mbox{-}\hspace{-0.04cm}\nargmin g = \big\{x\in \dom g~\big|~ g(x) \leq \inf g + \delta\big\},
	\]
	where $\inf g = \inf\{g(x)~|~x\in X\}$ and $\dom g = \{x\in X~|~g(x)<\infty\}$. Convergence $x^\nu$ to $x$ along a subsequence indexed by $N\subset\nats$ is specified by $x^\nu\Nto x$. The index set $N$ is then taken from the collection of subsequences of $\nats$ denoted by $\cN_\infty^\grill$.
	
	\begin{definition}\label{def:epiconv} {\rm (epi-convergence).}
		For $g^\nu,g:X\to \Reals$, we say that $g^\nu$ epi-converges to $g$ as $\nu \to \infty$, written $g^\nu\eto g$, when the following hold at every $x\in X$:
		\begin{enumerate}[(a)]
			\item For all $x^\nu\to x$, one has $\nliminf g^\nu(x^\nu) \geq g(x)$.
			
			\item There exists $x^\nu\to x$ such that $\nlimsup g^\nu(x^\nu) \leq g(x)$.
			
		\end{enumerate}
	\end{definition}
	
	The main result of this section furnishes an optimality gap for a control obtained by solving the approximating problems \eqref{eqn:approx}.
	
	\begin{theorem}\label{thm:minval} {\rm (optimality gap).} For a metric space $(X,d_X)$, suppose that $f:X\to \Reals$, $X^n\subset X$, and $T_n:\reals^n\to X^n$ is a surjective mapping. Let $f_n, f_n^\nu:\reals^n\to \Reals$, with $f_n$ defined by \eqref{eqn:fn}.  Suppose that the following hold for some $\bar n\in \nats$:
		\begin{enumerate}[(a)]
			
			\item $f + \iota_{X^n} \eto f$ as $n\to \infty$.
			
			\item $\forall n\geq \bar n$, $f_n^\nu\eto f_n$ as $\nu\to \infty$.
			
			\item $\inf f > -\infty$ and, $\forall n\geq \bar n$, $\inf f_n < \infty$.
		\end{enumerate}
		Then, for any $\epsilon \in (0,\infty)$, there exists $n_\epsilon\geq \bar n$ such that when $n\geq n_\epsilon$, $\delta^\nu \to \delta \in [0,\infty)$ and $\{\bar x_n^\nu\in \delta^\nu\mbox{-}\nargmin f_n^\nu, \nu\in\nats\}$ has a cluster point $\bar x_n$, one obtains
		\[
		f\big( T_n(\bar x_n)\big) \leq \inf f + \epsilon + \delta.
		\]
	\end{theorem}
	\state Proof. Suppose that $\inf f$ is finite and let $\gamma\in (0,\infty)$. Then, there is $x\in X$ such that $f(x) \leq \inf f + \gamma$. In view of (a) and Definition \ref{def:epiconv}(b), there is $x^n\to x$ such that $\nlimsup_n (f(x^n) + \iota_{X^n}(x^n)) \leq f(x)$. Stringing these two inequalities together, we obtain that
	\[
	\nlimsup_n \big(\inf (f + \iota_{X^n})\big)\leq \nlimsup_n \big(f(x^n) + \iota_{X^n}(x^n)\big) \leq f(x) \leq \inf f + \gamma.
	\]
	Since $\gamma$ is arbitrary, we conclude that
	\begin{equation}\label{eqn:proofOptgaplimsup}
		\nlimsup_n \big(\inf (f+\iota_{X^n})\big) \leq \inf f
	\end{equation}
	when $\inf f$ is finite. The same holds trivially when $\inf f = \infty$. The possibility $\inf f = -\infty$ is ruled out by (c).

	Let $\epsilon \in (0,\infty)$. By (c) and \eqref{eqn:proofOptgaplimsup}, $\inf f>-\infty$ and there is $n_\epsilon \geq \bar n$ such that
	\begin{equation}\label{eqn:infineq}
		\inf (f+\iota_{X^n}) \leq \inf f + \epsilon ~~~\forall n\geq n_\epsilon.
	\end{equation}
	
	Fix $n \geq n_\epsilon$. First, we establish that
	\begin{equation}\label{eqn:proofOptgaplimsup2}
		\nlimsup_\nu (\inf f_n^\nu) \leq \inf f_n.
	\end{equation}
	An argument parallel to the one carried out to reach \eqref{eqn:proofOptgaplimsup} confirms the claim when $\inf f_n$ is finite. Since $\inf f_n = \infty$ is ruled out by (c), it only remains to examine the case with $\inf f_n = -\infty$. Let $\gamma \in (0,\infty)$. Then, there is $x_n\in \reals^n$ such that $f_n(x_n) \leq -\gamma$ and also, by Definition \ref{def:epiconv}(b), there are $x_n^\nu\to x_n$ such that $\nlimsup_\nu f_n^\nu(x_n^\nu) \leq f_n(x_n)$. Thus,
	\[
	\nlimsup_\nu (\inf f_n^\nu) \leq \nlimsup_\nu f_n^\nu(x_n^\nu) \leq f_n(x_n) \leq -\gamma.
	\]
	Since $\gamma$ is arbitrary, \eqref{eqn:proofOptgaplimsup2} holds in this case as well.
	
	Second, let $N\in \cN_\infty^\grill$ be the subsequence corresponding to the cluster point $\bar x_n$, i.e., $\bar x_n^\nu \Nto \bar x_n$. Then, \eqref{eqn:proofOptgaplimsup2} and the construction of $\{\bar x_n^\nu, \nu\in\nats\}$ ensure that
	\begin{equation*}
		\nlimsup_{\nu\in N} f_n^{\nu}(\bar x_n^{\nu}) \leq \nlimsup_{\nu\in N} (\inf f_n^{\nu} + \delta^{\nu}) \leq \inf f_n + \delta.
	\end{equation*}
	By Definition \ref{def:epiconv}(a), this implies that
	\begin{equation*}
		f_n(\bar x_n) \leq \nliminf_{\nu\in N} f_n^{\nu}(\bar x_n^{\nu}) \leq \nlimsup_{\nu\in N} f_n^{\nu}(\bar x_n^{\nu})\leq \inf f_n + \delta.
	\end{equation*}
	Since $\inf f_n< \infty$ by (c), $\bar x_n\in \dom f_n$. Thus, $\bar x_n \in \delta\mbox{-}\hspace{-0.04cm}\nargmin f_n$. We then obtain
	\[
	f\big(T_n(\bar x_n)\big) = f_n(\bar x_n) \leq \inf f_n + \delta = \inf (f + \iota_{X^n}) + \delta \leq \inf f + \epsilon + \delta,
	\]
	where the second equality follows by \eqref{eqn:equivInf} and the second inequality by \eqref{eqn:infineq}.\eop

	Theorem \ref{thm:minval} provides an optimality gap, relative to the infinite-dimensional problem \eqref{eqn:actualpr}, for any cluster point constructed by solving a sequence of approximating finite-dimensional problems. The gap is bounded by the sum of two terms: $\epsilon$ and $\delta$. While $\epsilon$ can be selected arbitrarily near zero, this typically results in a large $n_\epsilon$ and thus the need for minimizing a high-dimensional function $f_n^\nu$ with a tolerance eventually near $\delta$. The theorem applies for any $\delta$ so the approximating problems do not need to be solved to a high tolerance. The tolerance might not even be fully known. Still, the theorem provides a guarantee that better solutions of the approximating problems (i.e., lower $\delta^\nu$) translate into better solutions of the infinite-dimensional problem \eqref{eqn:actualpr}, at least in the sense of an upper bound.
	
	The theorem does not require convexity, smoothness or continuity of the functions $f,f_n$, and $f_n^\nu$. However, the functions $f$ and $f_n$ are lower semicontinuous (lsc) by virtue of being epi-limits. We note that $\{\bar x_n^\nu, \nu\in \nats\}$ is a sequence in $\reals^n$ and thus has a cluster point as long as the sequence is bounded.
	
	In general, $\epsilon$ can be thought of as a {\em discretization error} due to the restriction from $X$ to $X^n$ and $\delta^\nu$ (and its limiting value $\delta$) is an {\em optimization error} caused by our incomplete minimization of $f_n^\nu$. There is no {\em approximation error} in the theorem, something one might expect from replacing $f_n$ in \eqref{eqn:discrF} by $f_n^\nu$ in  \eqref{eqn:approx}. This is avoided because, as $\nu\to \infty$, such errors vanish by assumption (b) in Theorem \ref{thm:minval}.

\section{Performance Bounds}\label{sec:performance-bounds}
	
We now return to PDE-constrained optimization under uncertainty (OUU) and the notation in Section 2. Following the pattern of Section 3 with \eqref{eqn:actualpr} being replaced by the problems \eqref{eqn:discr},  \eqref{eqn:discrF}, and \eqref{eqn:approx}, we start by defining the {\em control-discretized problems}
\begin{equation}\label{eqn:discrpde}
		\Big\{\nnmin_{z\in Z} ~\phi(z) + \iota_{Z^n}(z), ~~n\in \nats\Big\},
	\end{equation}
	where $\phi:Z\to \Reals$ is the objective function in the actual problem \eqref{eqn:actualpde0} and $Z^n\subset Z$. We assume that elements in $Z^n$ can be represented by $n$ real-valued parameters and that
	there is a surjective mapping $T_n:\reals^n\to  Z^n$. Through this correspondence, we obtain the {\em equivalent control-discretized problems}
	\begin{equation}\label{eqn:discrpdeF}
		\Big\{\nnmin_{z_n\in \reals^{n}} ~\phi_n(z_n), ~~n\in\nats\Big\},
	\end{equation}
	where $\phi_n:\reals^{n}\to \Reals$ is given by $\phi_n(z_n) = \phi( T_n(z_n))$. Leveraging the definition of $\phi$ in \eqref{eqn:actualpde0}, the problems in \eqref{eqn:discrpdeF} take the form
	\begin{equation}\label{eqn:discrpdeFexpl}
		\nnmin_{z_n\in \reals^{n}} ~\phi_n(z_n) = \iota_{A_n}(z_n) + f_n(z_n) + h\Big(\Ex\big[G_n(\bfxi,z_n)\big]\Big),
	\end{equation}
	where $A_n = \{z_n\in \reals^n~|~T_n(z_n) \in A\}$, $f_n(z_n) = f(T_n(z_n))$, and $G_n(\xi,z_n) = G(\xi, T_n(z_n))$ so that $f_n:\reals^n\to \reals$ and $G_n:\Xi\times \reals^n\to \reals^{m}$.
	
\begin{assumption}{\rm (continuity).}\label{ass:rlsc}
For a.e.\ $\xi\in \Xi$ and $i=1, \dots, m,$ $s(\xi,\cdot\,):Z\to U$, $g_i$, and $T_n$ are continuous\footnote{Throughout, a product space is equipped with the product norm.}. For every $z_n\in\reals^n$, $G_n(\cdot\,,z_n):\Xi\to \reals^m$ is measurable.
	\end{assumption}
	
	Under Assumption \ref{ass:rlsc}, $\Ex[G_n(\bfxi,z_n)]$ is well defined for each $z_n\in\reals^n$ under the usual convention of setting an expectation to infinity if both the positive part and the negative part integrate to infinity \cite[Section 3.B]{primer}.
	
	Next, we turn to approximations of the equivalent control-discretized problems. For $\xi_1, \dots, \xi_\nu\in \Xi$, the {\em approximating problems} take the form
	\begin{equation}\label{eqn:discrpdeFexplapprox}
		\nnmin_{z_n\in \reals^{n}} ~ \phi_n^\nu(z_n) = \iota_{A_n}(z_n) + f_n(z_n) + h^\nu\Bigg(\frac{1}{\nu}\sum_{j=1}^\nu G_n^\nu(\xi_j,z_n)\Bigg),
	\end{equation}
where a sample average  replaces the expectation, $h^\nu:\reals^{m}\to \Reals$ approximates the monitoring function, and $G_n^\nu:\Xi\times\reals^n\to \reals^{m}$ approximates $G_n$ and is given by
	\begin{equation}\label{eq:G-n}
		G_n^\nu(\xi,z_n) = G^\nu\big(\xi, T_n(z_n)\big).
	\end{equation}
Here, $G^\nu:\Xi\times Z\to \reals^m$ has
	\begin{equation}\label{eq:G-nu}
			G^\nu(\xi,z) = \Big( g_1^\nu\big( s^{\nu}(\xi,z), z\big), \dots, g_m^\nu\big( s^{\nu}(\xi,z), z\big)\Big),
	\end{equation}
	where $g_i^\nu:U\times Z\to \reals$ approximates $g_i$ and $s^\nu:\Xi\times Z\to U$ approximates  $s:\Xi\times Z\to U$; see Section \ref{sec:applications}. We consider the following scheme to solve \eqref{eqn:actualpde0}.

	\bigskip
	
	\state Approximation Algorithm for PDE-Constrained OUU.
	
	\begin{description}
		
		\item[Data.] ~~$n\in \nats$ and $\{\xi_j\in \Xi, j\in\nats\}$.

		\item[Step 0.] Set $\nu=1$.

		\item[Step 1.] Apply an algorithm to \eqref{eqn:discrpdeFexplapprox}, obtain $\bar z_n^\nu\in\reals^n$ and record $\phi_n^\nu(\bar z_n^\nu)$.
		
		\item[Step 2.] Replace $\nu$ by $\nu+1$ and go to Step 1.
	\end{description}
	
	\smallskip

	Presumably, the approximating problems \eqref{eqn:discrpdeFexplapprox} are tractable for existing (nonlinear programming) algorithms, but there is no requirement that $\bar z_n^\nu$ must be globally or locally optimal or even stationary for $\phi_n^\nu$. A main feature of the subsequent analysis, already alluded to in Section 3, is that the resulting performance guarantee holds even though Step 1 is not carried out ``perfectly.'' In particular, if $\{\bar z_n^\nu, \nu\in \nats\}$ has a cluster point $\bar z_n$, which indeed would be the case when the sequence is bounded, then the recorded values $\{\phi_n^\nu(\bar z_n^\nu), \nu\in \nats\}$ lead to a performance guarantee for $\bar z_n$ as measured by the actual objective function $\phi$; see Theorem \ref{thm:perfguar} below.
	
	Numerous implementation issues emerge including rules for refining approximation levels and stopping criteria in Step 1. It would typically be inefficient to apply much effort in Step 1 for {\em each} $\nu$. In practice, one might leapfrog over most $\nu$ and invest computing resources towards, say, $\nu = 1000, 2000, 3000$, etc. Each of the expectations might also demand different sample sizes. They are all set to $\nu$ in this paper for notational convenience, but the results extend to other schemes trivially. While $n$ is kept fixed in the algorithm, an implementation might also involve a gradual increase of that parameter.

\subsection{Intermediate Results}\label{sec:intermediate}
	
	This subsection presents assumptions and furnishes several technical results. In the following, inequalities between vectors are assumed to apply componentwise. The Euclidean balls in $\reals^n$ are written as
	\[
	\ball(\bar z_n,\rho) = \big\{z_n\in\reals^n~\big|~\|z_n - \bar z_n\|_2\leq \rho\big\}.
	\]
	
\begin{assumption}{\rm (locally bounded quantities of interest).}\label{ass:Gbound} For each $\bar z_n\in A_n$, there are $\rho\in (0,\infty)$ and $P$-integrable $c:\Xi\to [0,\infty)$ such that
		\[
		\big\|G_n(\xi, z_n)\big\|_\infty \leq c(\xi) ~~~\forall z_n\in \ball(\bar z_n,\rho)\cap A_n, ~~\mbox{a.e. } \xi \in \Xi.
		\]
\end{assumption}

	\begin{lemma}\label{lem:epiLLN} Let $n\in \nats$ be fixed. Suppose that $A_n$ is closed, Assumptions \ref{ass:rlsc} and \ref{ass:Gbound} hold for $n$ and the sample $\{\bfxi_1, \bfxi_2, \dots\}$ is iid as $\bfxi$. Then, with probability one,
		\[
		\nlim_\nu \frac{1}{\nu}\sum_{j=1}^\nu G_n(\bfxi_j, z_n^\nu)  = \Ex\big[G_n(\bfxi, z_n)\big] \in \reals^m~~\mbox{ whenever } z_n^\nu \in A_n \to z_n.
		\]
	\end{lemma}
	\state Proof. Let $\psi_i:\Xi\times\reals^n\to \reals$ be given by $\psi_i(\xi,z_n) = g_i(s(\xi,T_n(z_n)),T_n(z_n))$. For a.e. $\xi\in \Xi$, $\psi_i(\xi, \cdot\,)$ is continuous by Assumption \ref{ass:rlsc}. For all $z_n\in \reals^n$, $\psi_i(\cdot, z_n)$ is measurable. Thus, $\psi_i$ is a Caratheodory function and certainly random lsc; see \cite[Ex. 8.51(d)]{primer}. Moreover, $(\xi,z_n)\mapsto \psi_i(\xi,z_n) + \iota_{A_n}(z_n)$ is random lsc by \cite[Ex. 8.51(c)]{primer} because $A_n$ is closed.
With probability one, the epigraphical law of large numbers \cite[Thm. 8.56]{primer} applied to each $\psi_i$ yields
\[
\nliminf_\nu \frac{1}{\nu}\sum_{j=1}^\nu G_n(\bfxi_j, z_n^\nu)  \geq \Ex\big[G_n(\bfxi, z_n)\big] ~~\mbox{ whenever } z_n^\nu \in A_n \to z_n,
\]
with no component of $\Ex[G_n(\bfxi, z_n)]$ being $-\infty$. Repeating this argument with $-\psi_i$ in place of $\psi_i$, we obtain the conclusion.\eop
	
\begin{assumption}{\rm (approximation of performance functions).}\label{ass:smooth}
There is a sequence $\{\epsilon^\nu \in [0,\infty), ~\nu\in\nats\}$ with $\epsilon^\nu \to 0$ such that
\[
\big| g_i^\nu(u,z) - g_i(u,z) \big| \leq \epsilon^\nu~~~~\forall u\in U, ~z\in Z^n\cap A, ~i=1, \dots, m.
\]
\end{assumption}

The assumption requires a uniform error bound for the approximating performance functions, which is satisfied, for example by the smax-approximation in Subsection \ref{subsec:approx}. The uniformity can be relaxed, with slight adjustments to Assumption \ref{ass:upper} below. We omit these details due to the resulting notational complexity. Anyhow, they are not needed for the application in Section \ref{sec:applications}.

\begin{assumption}{\rm (solution properties and approximations).}\label{ass:upper}
There exist $\eta,\kappa\in (0,\infty)$, $\lambda,\mu:\Xi\to [0,\infty)$, $\rho,\pi:A_n\to (0,\infty)$, and $\{\Delta_0,\Delta_\nu:A_n\to [0,\infty),~\nu\in\nats\}$ satisfying the following properties:

\hspace{0.01cm}

\begin{enumerate}[(a)]

\item $\Ex\big[\lambda^\eta(\bfxi)\mu^\kappa(\bfxi)\big]$ and $\,\Ex\big[\mu^\kappa(\bfxi)\big]$ are finite.

\hspace{0.01cm}

\item $\forall z_n\in A_n$, $z^n = T_n(z_n)$, $u,u'\in U$, $\nu\in\nats$, $i= 1, \dots, m$, and a.e. $\xi\in\Xi$, one has
\begin{align*}
\big| g_i(u,z^n) - g_i(u',z^n) \big|  & \leq \pi(z_n) \big(\|u\|_U^\eta + \|u'\|_U^\eta + 1\big) \|u - u'\|_U^\kappa\\
\max\Big\{ \big\|s(\xi,z^n)\big\|_U, ~~ \big\|s^\nu(\xi,z^n)\big\|_U\Big\} & \leq \lambda(\xi) \Delta_0(z_n)\\
\big\| s^\nu(\xi, z^n ) - s(\xi, z^n ) \big\|_U & \leq \mu(\xi)\Delta_\nu(z_n).
\end{align*}

\item $\forall\bar z_n\in A_n$, one has
\begin{align*}
\sup_{z_n\in \ball(\bar z_n,\rho(\bar z_n))\cap A_n} & \pi(z_n) < \infty, ~~~~~~~\sup_{z_n\in \ball(\bar z_n,\rho(\bar z_n))\cap A_n} & \Delta_0(z_n) < \infty,\\
\sup_{z_n\in \ball(\bar z_n,\rho(\bar z_n))\cap A_n} & \Delta_\nu(z_n) \to 0, ~~\mbox{ as } \nu\to \infty.
\end{align*}
\end{enumerate}
\end{assumption}

Part (b) imposes bounds on approximation errors of the solutions as well as a relatively mild H\"{o}lder-type condition on the performance functions, which can be made more general depending on specific applications. Section \ref{sec:applications} furnishes a specific example of when these assumptions hold.
	
\begin{lemma}\label{lem:nuerrorlemma} Let $n\in \nats$ be fixed. Suppose that $A_n$ is closed, Assumptions \ref{ass:smooth} and \ref{ass:upper} hold for $n$ and the sample $\{\bfxi_1, \bfxi_2, \dots\}$ is iid at $\bfxi$. Then, with probability one,
		\[
		\frac{1}{\nu}\sum_{j=1}^{\nu} \big\|G_n^\nu(\bfxi_j,z_n^\nu) - G_n(\bfxi_j,z_n^\nu) \big\|_\infty \to 0 ~~\mbox{ whenever}~~z_n^\nu \in A_n \to z_n.
		\]
	\end{lemma}
\state Proof. Let $\eta,\kappa$, $\lambda,\mu$, $\rho$, $\pi$, $\Delta_0$, $\Delta_\nu$, and $\epsilon^\nu$ be as specified by Assumptions \ref{ass:smooth} and \ref{ass:upper}. First, we fix $z_n \in A_n$, $z^n = T_n(z_n)$, $\xi \in \Xi_1$, $\nu\in \nats$, and $i\in \{1, \dots, m\}$, where $P(\Xi_1) = 1$. Trivially,
\begin{align*}
  & \Big| g_i^\nu\big(s^\nu(\xi,z^n),z^n\big) - g_i\big(s(\xi,z^n),z^n\big)\Big|\\
   & \leq \Big| g_i\big(s^\nu(\xi,z^n),z^n\big) - g_i\big(s(\xi,z^n),z^n\big) \Big| + \Big| g_i^\nu\big(s^\nu(\xi,z^n),z^n\big) - g_i\big(s^\nu(\xi,z^n),z^n\big) \Big|.
\end{align*}
Assumption \ref{ass:upper} addresses the first term so that
\begin{align*}
  & \Big| g_i\big( s^\nu(\xi,z^n), z^n\big) - g_i\big( s(\xi,z^n), z^n\big) \Big| \\
  & \leq  \pi(z_n) \Big(\big\|s^\nu(\xi,z^n)\big\|_U^\eta + \big\|s(\xi,z^n)\big\|_U^\eta + 1\Big) \big\|s^\nu(\xi,z^n) - s(\xi,z^n)\big\|_U^\kappa\\
  & \leq  \pi(z_n) \Big(2 \lambda^\eta(\xi) \Delta_0^\eta(z_n) +1\big)\mu^\kappa(\xi)\Delta_\nu^\kappa(z_n).
\end{align*}
Also bringing in Assumption \ref{ass:smooth}, we obtain for $\{\xi_j\in \Xi_1, j\in \nats\}$ that
\begin{align*}
\frac{1}{\nu}\sum_{j=1}^\nu \big\| G_n^\nu(\xi_j,z_n) - G_n(\xi_j,z_n)\big\|_\infty & \leq  \pi(z_n) 2 \Delta_0^\eta(z_n) \Delta_\nu^\kappa(z_n) \frac{1}{\nu}\sum_{j=1}^\nu \lambda^\eta(\xi_j)  \mu^\kappa(\xi_j)\\
&~~+  \pi(z_n) \Delta_\nu^\kappa(z_n) \frac{1}{\nu}\sum_{j=1}^\nu \mu^\kappa(\xi_j) + \epsilon^\nu.
\end{align*}

Second, Assumption \ref{ass:upper}(a) and the iid sampling ensure that
\[
\frac{1}{\nu}\sum_{j=1}^\nu \lambda^\eta(\bfxi_j)\mu^\kappa(\bfxi_j) \to \Ex\big[\lambda^\eta(\bfxi)\mu^\kappa(\bfxi)\big]\in\reals, ~~~~ \frac{1}{\nu}\sum_{j=1}^\nu \mu^\kappa(\bfxi_j) \to \Ex\big[\mu^\kappa(\bfxi)\big]\in\reals
\]
with probability one. Let $\xi_1, \xi_2, \dots$ be a sequence in $\Xi_1$ for which this convergence holds and let $z_n^\nu \in A_n \to \bar z_n$. Then, there exists $\bar\nu$ such that $z_n^\nu\in \ball( \bar z_n,\rho(\bar z_n))$ for all $\nu\geq\bar\nu$. Let $B = \ball( \bar z_n,\rho(\bar z_n)) \cap A_n$, $\bar\pi = \sup_{z_n \in B} \pi(z_n)$, $\bar\Delta = \sup_{z_n \in B} \Delta_0(z_n)$, and $\delta_\nu = \sup_{z_n \in B} \Delta_\nu(z_n)$.  Thus, for all $\nu\geq \bar\nu$,
\[
\frac{1}{\nu}\sum_{j=1}^\nu \big\| G_n^\nu(\xi_j,z_n^\nu) - G_n(\xi_j,z_n^\nu)\big\|_\infty \leq  \frac{ \bar\pi 2 \bar\Delta^\eta \delta_\nu^\kappa}{\nu}\sum_{j=1}^\nu \lambda^\eta(\xi_j)  \mu^\kappa(\xi_j) +  \frac{\bar\pi \delta_\nu^\kappa}{\nu}\sum_{j=1}^\nu \mu^\kappa(\xi_j) + \epsilon^\nu.
\]
Since $\bar\pi,\bar\Delta<\infty$ and $\delta_\nu,\epsilon^\nu\to 0$, the conclusion follows.\eop

\begin{assumption}{\rm (locally bounded quantities of interest on $Z$).}\label{ass:GboundInf}
		For each $\bar z\in A$, $G(\cdot, \bar z)$ is measurable and there are $\rho\in (0,\infty)$ and $P$-integrable $c:\Xi\to [0,\infty)$ such that
		\[
		\big\|G(\xi, z)\big\|_\infty \leq c(\xi) ~~~\mbox{ for } ~z\in A, ~\|z - \bar z\|_Z \leq \rho, ~~\mbox{a.e. }\xi \in \Xi.
		\]
\end{assumption}
	
This assumption resembles Assumption \ref{ass:Gbound}, but is stated separately to avoid intricacies about the relation between $Z$ and $\reals^n$.

\begin{lemma}{\rm (continuity of expectation function).}\label{lem:llnInfDim}
If $A$ is closed and Assumptions \ref{ass:rlsc} and \ref{ass:GboundInf} hold, then $\Ex[G(\bfxi,z^n)] \to \Ex[G(\bfxi,z)]\in\reals^m$ whenever $z^n\in A\to z$.
\end{lemma}
\state Proof. Under the stated assumptions, $G(\xi,\cdot\,)$ is continuous for a.e.\ $\xi$. Thus, the fact follows from a standard application of the dominated convergence theorem.\eop

\subsection{Main Results}
	
We establish two performance guarantees related to the Approximation Algorithm for PDE-Constrained OUU. The first one furnishes an upper bound on the objective function value in the actual problem \eqref{eqn:actualpde0} for any cluster point produced by the algorithm. The second one specifies an optimality gap.
	
For metric space $X$, a function $g:X\to \Reals$ is lsc relative to $C\subset X$ if $\nliminf g(x^\nu) \geq g(x)$ whenever $x^\nu\in C\to x\in C$.

\begin{theorem}{\rm (upper bound in PDE-constrained OUU).}\label{thm:perfguar}
For fixed $n\in\nats$, suppose that $A_n$ is closed, $f_n$ is lsc relative to $A_n$, the sample $\{\bfxi_1, \bfxi_2, \dots\}$ is iid as $\bfxi$, Assumptions \ref{ass:rlsc}, \ref{ass:Gbound}, \ref{ass:smooth}, and \ref{ass:upper} hold for $n$ and
		\[
		\nliminf h^\nu(w^\nu) \geq h(w)>-\infty~~~\mbox{ whenever } w^\nu \to w\in\reals^m.
		\]
With probability one, if $\{\bar z_n^\nu \in \reals^{n}, \nu\in\nats\}$ converges to $\bar z_n$ along a subsequence $N\in \cN_\infty^\grill$ as $\nu\to \infty$, then
		\[
		\phi\big( T_n(\bar z_n)\big) \leq \nliminf_{\nu\in N} \phi_n^\nu(\bar z_n^\nu).
		\]
\end{theorem}
\state Proof. By Lemma \ref{lem:epiLLN} and Lemma \ref{lem:nuerrorlemma}, with probability one, $z^\nu_n \in A_n\to z_n$ implies
\begin{equation}\label{eqn:convprop}
\begin{split}
&\frac{1}{\nu} \sum_{j=1}^{\nu} G_n(\bfxi_j,z_n^\nu) \to \Ex\big[G_n(\bfxi,z_n)\big]\in\reals^m\\
&\frac{1}{\nu} \sum_{j=1}^\nu \big\|G_n^\nu(\bfxi_j,z_n^\nu) - G_n(\bfxi_j,z_n^\nu)\big\|_\infty \to 0.
\end{split}
\end{equation}
Let $\{\xi_1, \xi_2, \dots\}$ be an event for which this occurs. Then,
\[
\frac{1}{\nu} \sum_{j=1}^{\nu} G_n^\nu(\xi_j,z_n^\nu) = \frac{1}{\nu} \sum_{j=1}^\nu \big(G_n^\nu(\xi_j, z_n^\nu) - G_n(\xi_j, z_n^\nu)\big) + \frac{1}{\nu} \sum_{j=1}^{\nu} G_n(\xi_j, z_n^\nu) \to \Ex\big[G_n(\bfxi, z_n)\big]
\]
whenever $z^\nu_n \in A_n\to z_n$.

Suppose that $\bar z_n^\nu\Nto \bar z_n$ for $N\in \cN_\infty^\grill$. We aim to show that
\begin{equation}\label{eqn:liminfcond}
	\nliminf_{\nu\in N} \phi_n^\nu(\bar z_n^\nu) \geq \phi_n(\bar z_n)
\end{equation}
because then the conclusion follows from the fact that $\phi_n(\bar z_n) = \phi(T_n(\bar z_n))$.

If $\bar z_n \not\in A_n$, then $\bar z^\nu_n \not\in A_n$ for sufficiently large $\nu \in N$ because $A_n$ is closed. Thus, \eqref{eqn:liminfcond} holds trivially in this case and we concentrate on the case with $\bar z_n \in A_n$. Without loss of generality, we assume that $\bar z_n^\nu \in A_n$ for all $\nu\in N$ because any $\bar z_n^\nu \not\in A_n$ produces $\phi_n^\nu(\bar z_n^\nu) = \infty$.  These facts and the assumption on the monitoring functions imply that
\[
\nliminf_{\nu \in N} h^\nu\Bigg( \frac{1}{\nu} \sum_{j=1}^{\nu} G_n^\nu(\xi_j,\bar z_n^\nu) \Bigg) \geq h\Big(\Ex\big[G_n(\bfxi,\bar z_n)\big]\Big)>-\infty.
\]
Since $f_n$ is lsc relative to $A_n$, one has
\[
\nliminf_{\nu\in N} f_n(\bar z_n^\nu) \geq f_n(\bar z_n).
\]
All the three terms defining $\phi_n^\nu(\bar z_n^\nu)$ are bounded from below. Thus, it follows that
\begin{align*}
		&\nliminf_{\nu\in N} \phi_n^\nu(\bar z_n^\nu)\\
		& \geq \nliminf_{\nu\in N} \iota_{A_n}(\bar z_n^\nu) + \nliminf_{\nu\in N} f_n(\bar z_n^\nu) + \nliminf_{\nu \in N} h^\nu\Bigg( \frac{1}{\nu} \sum_{j=1}^{\nu} G_n^\nu(\xi_j,\bar z_n^\nu) \Bigg)\\
		& \geq \iota_{A_n}(\bar z_n) + f_n(\bar z_n) + h\Big(\Ex\big[G_n(\bfxi,\bar z_n)\big]\Big) = \phi_n(\bar z_n),
\end{align*}
which completes the proof.\eop

Theorem \ref{thm:perfguar} implies that if the algorithm constructs a sequence $\{\bar z_n^\nu, \nu\in\nats\}$ with a cluster point, then that point---a finite-dimensional vector---corresponds to a point in the control space $Z$, which is as good as we would expect from the recorded values $\{\phi_n^\nu(\bar z_n^\nu), \nu\in \nats\}$.  The theorem holds regardless of $n$, but a large $n$ would typically be associated with a large set $Z^n$ and thus better chances to obtain low values $\phi_n^\nu(\bar z_n^\nu)$.

The assumption on $h^\nu$ and $h$ is satisfied by models where inequality constraints are replaced by penalties. For example, $h(w) = \iota_{(-\infty,0]^m}(w)$ is approximated by $h^\nu(w)$ $=$ $\theta^\nu \sum_{i=1}^m (\max\{0,w_i\})^2$, where $\theta^\nu\to \infty$.
	
	There is no need to solve \eqref{eqn:discrpdeFexplapprox} to local or global optimality. One would simply attempt to bring the objective function value down as $\nu\to \infty$. In practice, one might have a goal of obtaining a control $z\in Z$ such that $\phi(z) \leq \alpha$. The theorem provides a way of certifying this: pick a reasonably large $n$, apply the algorithm with a stopping criterion in Step 1 of $\phi_n^\nu(\bar z_n^\nu)\leq \alpha$. Any cluster point $\bar z_n$ of the constructed sequence produces a control $\bar z = T_n(\bar z_n)$ which then is good enough. Since the theorem does not leverage any assumptions on the discretization $Z^n$ of $Z$, it becomes impossible to certify a priori whether a good solution of this kind is achievable; one needs to wait for the computations and see what values come out.

	Next, we introduce additional assumptions and these allow us to claim that for sufficiently large $n$, the algorithm obtains a decision that is arbitrarily good relative to the actual problem.
	
	In the following, the $i$th component of a vector $v$ is sometimes indicated by $(v)_i$.

\begin{assumption}{\rm (constraint qualification).}\label{ass:CQ}
For fixed $n$ and all points $z_n\in A_n$ satisfying $h(\Ex[G_n(\bfxi,z_n)])$ $<$ $\infty$, there exists $z^\nu_n \in A_n\to z_n$ as $\nu\to \infty$ such that
\[			
\Big(\Ex\big[G_n(\bfxi,z^\nu_n)\big]\Big)_i < \Big(\Ex\big[G_n(\bfxi,z_n)\big]\Big)_i, ~~i=1, \dots, m, ~\nu\in\nats.
\]

For all $z\in A$, with $h(\Ex[G(\bfxi,z)])$ $<$ $\infty$, there exists $z^n \in A \cap Z^n\to z$ as $n\to \infty$ such that $\Ex[G(\bfxi,z^n)] \leq \Ex[G(\bfxi,z)]$ for all $n\in\nats$.
\end{assumption}

\begin{theorem}{\rm (optimality gap in PDE-constrained OUU).}\label{thm:optgap}
For fixed $\bar n\in\nats$, suppose that $\inf \phi > -\infty$, $A$ is closed, $f$ is continuous relative to $A$, the sample $\{\bfxi_1, \bfxi_2, \dots\}$ is iid as $\bfxi$, and, for all $n\geq \bar n$, $\inf \phi_n < \infty$, and suppose that Assumptions \ref{ass:rlsc}, \ref{ass:Gbound}, \ref{ass:smooth}, \ref{ass:upper}, and \ref{ass:GboundInf} are satisfied. Moreover, either (a) or (b) holds:
\begin{enumerate}[(a)]
			
\item $h$ is continuous and $h^\nu(w^\nu) \to h(w) > -\infty$  whenever $w^\nu \to w\in\reals^m$.
			
$\forall z\in A$, there exists $z^n \in  A \cap Z^n\to z$.

\item $\nliminf_\nu h^\nu(w^\nu) \geq h(w)>-\infty$ whenever $w^\nu \to w\in\reals^m$.
			
$h^\nu(w) \leq h(w)$ for all $w\in\reals^m$, $\nu\in\nats$.

$h$ is lsc and $h(w) \leq h(w')$ for all $w\leq w'$.

Assumption \ref{ass:CQ} is satisfied for all $n\geq \bar n$.
			
\end{enumerate}
With probability one, the following holds: for any $\epsilon \in (0,\infty)$, there exists $n_\epsilon\geq \bar n$ such that if $n\geq n_\epsilon$, $\delta^\nu \to \delta \in [0,\infty)$, and $\{\bar z_n^\nu\in \delta^\nu\mbox{-}\nargmin \phi_n^\nu, \nu\in\nats\}$ has a cluster point $\bar z_n$, then
\[
\phi\big( T_n(\bar z_n)\big) \leq \inf \phi + \epsilon + \delta.
\]
That is, the optimality gap for $\bar z_n$ is $\epsilon + \delta$.
\end{theorem}
\state Proof. Let $n\geq \bar n$. Since $A$ is closed and $T_n$ is continuous by Assumption \ref{ass:rlsc}, it follows that $A_n$ is closed. Likewise, the continuity of $f$ relative to $A$ implies that $f_n$ is continuous relative to $A_n$.

By Lemma \ref{lem:epiLLN} and Lemma \ref{lem:nuerrorlemma}, with probability one, any $z_n^\nu\in A_n\to z_n$ as $\nu\to \infty$ implies \eqref{eqn:convprop}. Let $\{\xi_1, \xi_2, \dots\}$ be an event for which this occurs. Since we consider a countable number of values of $n$, the probability-one set can be assumed to be independent of $n$.

By the argument in the proof of Theorem \ref{thm:perfguar}, we obtain that
	\[
	\nliminf_{\nu} \phi_n^\nu(z_n^\nu) \geq \phi_n(z_n) \mbox{ whenever } z_n^\nu \to z_n \mbox{ as } \nu\to \infty.
	\]
	We next show that for each $\hat z_n\in A_n$, there exists $\hat z_n^\nu\in A_n\to \hat z_n$ as $\nu\to \infty$ such that
	\[
	\nlimsup_{\nu} \phi_n^\nu(\hat z_n^\nu) \leq \phi_n(\hat z_n).
	\]
	This holds trivially if $\hat z_n \not\in A_n$. For $\hat z_n\in A_n$, we argue as follows.
	As in the proof of Theorem \ref{thm:perfguar}, $z_n^\nu \in A_n \to z_n$ as $\nu\to \infty$ implies that
	\begin{equation}\label{eqn:convsampl2}
	\frac{1}{\nu} \sum_{j=1}^{\nu} G_n^\nu(\xi_j,z_n^\nu) \to \Ex\big[G_n(\bfxi,z_n)\big], ~~ \mbox{ as } \nu\to \infty.
	\end{equation}
Now, if assumption (a) holds, then set $\hat z_n^\nu = \hat z_n$ for all $\nu$, which implies that
	\begin{align*}
		\nlimsup_{\nu} \phi_n^\nu(\hat z_n^\nu) & = f_n(\hat z_n) + \nlimsup_\nu h^\nu \Bigg( \frac{1}{\nu} \sum_{j=1}^\nu G_n^\nu(\xi_j,\hat z_n) \Bigg)\\
		& = f_n(\hat z_n) + h\Big(\Ex\big[G_n(\bfxi,\hat z_n)\big]\Big) = \phi_n(\hat z_n).
	\end{align*}
	Alternatively, if assumption (b) holds, then there are two cases. (i) If $h(\Ex[G_n(\bfxi,\hat z_n)])$ $=$ $\infty$, then we again set $\hat z_n^\nu = \hat z_n$ for all $\nu$ and obtain
	\[
	\nlimsup_{\nu} \phi_n^\nu(\hat z_n) \leq \phi_n(\hat z_n)
	\]
	because the right-hand side equals infinity. (ii)  If $h(\Ex[G_n(\bfxi,\hat z_n)])<\infty$, then, by the first part of Assumption \ref{ass:CQ}, there exist $\gamma^k>0$ and $\bar z_n^k \in A_n\to \hat z_n$ as $k\to \infty$ such that
	\[
	\Big(\Ex\big[G_n(\bfxi,\bar z_n^k)\big]\Big)_i + \gamma^k \leq \Big(\Ex\big[G_n(\bfxi,\hat z_n)\big]\Big)_i, ~~ i=1, \dots, m, ~~k\in\nats.
	\]
For each $k$, there exists by \eqref{eqn:convsampl2} $\nu_k$ such that for all $\nu\geq \nu_k$ one has
\[
	\Bigg(\frac{1}{\nu} \sum_{j=1}^{\nu} G_n^\nu(\xi_j,\bar z_n^k)\Bigg)_i \leq \Big(\Ex\big[G_n(\bfxi,\bar z_n^k)\big]\Big)_i + \gamma^k, ~~ i=1, \dots, m.
\]
These inequalities mean that, for each $\nu$, there is $k(\nu)$ such that
	\[
	\Bigg(\frac{1}{\nu} \sum_{j=1}^{\nu} G_n^\nu(\xi_j,\bar z_n^{k(\nu)})\Bigg)_i \leq
	\Big(\Ex\big[G_n(\bfxi,\hat z_n)\big]\Big)_i, ~~ i=1, \dots, m
	\]
and this can be done such that $k(\nu)\to \infty$ as $\nu\to \infty$. Set $\hat z_n^\nu = \bar z_n^{k(\nu)}$, which then tends to $\hat z_n$ as $\nu\to \infty$.
These facts and the assumption on $h^\nu$ and $h$ imply that
	\[
	\nlimsup_{\nu} h^\nu\Bigg( \frac{1}{\nu} \sum_{j=1}^{\nu} G_n^\nu(\xi_j,\hat z_n^\nu) \Bigg) \leq \nlimsup_{\nu} h\Bigg( \frac{1}{\nu} \sum_{j=1}^{\nu} G_n^\nu(\xi_j,\hat z_n^\nu) \Bigg)  \leq h\Big(\Ex\big[G_n(\bfxi,\hat z_n)\big]\Big)
	\]
and then also
	\begin{align*}
		\nlimsup_{\nu} \phi_n^\nu(\hat z_n^\nu) & \leq \nlimsup_\nu f_n(\hat z_n^\nu) + \nlimsup_\nu h^\nu \Bigg( \frac{1}{\nu} \sum_{j=1}^\nu G_n^\nu(\xi_j,\hat z_n^\nu) \Bigg)\\
		& \leq f_n(\hat z_n) + h\Big(\Ex\big[G_n(\bfxi,\hat z_n)\big]\Big) = \phi_n(\hat z_n).
	\end{align*}
Consequently, with probability one, $\phi_n^\nu \eto \phi_n$ as $\nu\to \infty$; see Definition \ref{def:epiconv}.
	
	Next, we confirm that $\phi + \iota_{Z^n} \eto \phi$ as $n\to \infty$. Let $z^n\to z$. If $z\not\in A$, then $z^n\not\in A$ for sufficiently large $n$ because $A$ is closed. Thus,
	\begin{equation}\label{eqn:liminfbdinfty}
		\nliminf_n \big(\phi(z^n) + \iota_{Z^n}(z^n)\big) \geq \phi(z)
	\end{equation}
	due to the fact that both sides equal infinity. If $z\in A$, then we argue as follows. If $z^n$ remains outside $A$ and/or $Z^n$, then the same inequality trivially holds because $\phi(z^n) + \iota_{Z^n}(z^n) = \infty$. Thus, we assume without loss of generality that $z^n\in A\cap Z^n$. Now,
	\begin{align*}
		\nliminf_n \big(\phi(z^n) + \iota_{Z^n}(z^n)\big) & = \nliminf_n \Big( f(z^n) + h\Big(\Ex\big[G(\bfxi,z^n)\big]\Big)\Big)\\
		& \geq  f(z) + \nliminf_n h\Big(\Ex\big[G(\bfxi,z^n)\big]\Big).
	\end{align*}
	In view of Lemma \ref{lem:llnInfDim}, $\Ex[G(\bfxi,z^n)] \to \Ex[G(\bfxi,z)]$. Under assumption (a), $h$ is continuous so that the liminf expression tends to  $h(\Ex[G(\bfxi,z)])$. Under assumption (b), $h$ is lsc and the same expression is bounded from below by  $h(\Ex[G(\bfxi,z)])$. In any case, we have confirmed that \eqref{eqn:liminfbdinfty} holds.
	
	Let $z\in Z$. We also need to construct $z^n\to z$ such that
	\begin{equation}\label{eqn:limsupcondMain}
	\nlimsup_n \big(\phi(z^n) + \iota_{Z^n}(z^n)\big) \leq \phi(z).
	\end{equation}
	We can assume that $z\in A$ because otherwise the right-hand side will be infinity. Under assumption (a), there is $z^n\in A \cap Z^n \to z$, which then implies that
	\begin{align*}
		\nlimsup_n \big(\phi(z^n) + \iota_{Z^n}(z^n)\big) & = \nlimsup_n \Big(f(z^n) + h\big(\Ex\big[G(\bfxi,z^n)\big]\big)\Big)\\
		& \leq f(z) + \nlimsup_n h\Big(\Ex\big[G(\bfxi,z^n)\big]\Big) = \phi(z),
	\end{align*}
	where we again leverage Lemma \ref{lem:llnInfDim}. Under assumption (b), we consider two cases. If $h(\Ex[G(\bfxi,z)])=\infty$, then $\phi(z) = \infty$ and \eqref{eqn:limsupcondMain} holds trivially. If  $h(\Ex[G(\bfxi,z)])<\infty$, then, by the second part of Assumption \ref{ass:CQ}, there is $z^n\in  A \cap Z^n\to z$ such that $\Ex[G(\bfxi,z^n)] \leq \Ex[G(\bfxi,z)]$. Since $h(w) \leq h(w')$ when $w\leq w'$, this implies that
	\begin{align*}
		\nlimsup_n \big(\phi(z^n) + \iota_{Z^n}(z^n)\big) & \leq f(z) + \nlimsup_n h\Big(\Ex\big[G(\bfxi,z^n)\big]\Big)\\
		& \leq f(z) + h\Big(\Ex\big[G(\bfxi,z)\big]\Big) = \phi(z).
	\end{align*}
We can conclude that $\phi + \iota_{Z^n} \eto \phi$; cf. Definition \ref{def:epiconv}.
	
We have satisfied the requirements of Theorem \ref{thm:minval} and the conclusion follows. \eop
	
Theorem \ref{thm:optgap} supplements Theorem \ref{thm:perfguar} by imposing conditions under which the Approximation Algorithm for PDE-Constrained OUU produces a solution with a specific optimality gap. The optimality gap consists of $\epsilon$, which can be made arbitrarily small at the expense of larger approximating problems in Step 1, and $\delta$, the tolerance invoked in Step 1.
	
The main additional assumptions relate to $h$ and its approximation $h^\nu$ as well as $Z^n$. Naturally, $Z^n$ needs to approximate $Z$ arbitrarily well as $n\to \infty$ and this suffices if $h$ is continuous; see (a) in the theorem. Trivially, $h$ is continuous when $h(w) = w\in\reals$ as in Example \ref{eExpectationReg} and when $h$ models finite weights and penalties applied to the various expectations. Continuity fails, however, when the monitoring function models inequality constraints by means of indicator functions. These situations are addressed via (b), where $Z^n$ in interplay with the expectation functions allows one to approach a feasible point along points that are {\em strictly} feasible; see Assumption \ref{ass:CQ}. This is a constraint qualification that resembles the Slater condition from convex optimization. In particular, (b) is tailored to situations when $h(w) = \iota_{(-\infty,0]^m}(w)$
% and its approximation is obtained through
with
 the penalization approximation $h^\nu(w) = \theta^\nu \sum_{i=1}^m (\max\{0,w_i\})^2$, where $\theta^\nu\to \infty$.

The main practical take away from the theorem is that one should not shy away from ``complicated'' formulations involving multiple expectation functions appearing in the objective function and in the constraints and even as part of compositions with other functions. Under relatively mild assumptions, the approximating optimization problems arising from discretization, sampling, penalization, smoothing, and other inaccuracies indeed produce solutions with optimality gap no more than $\epsilon$, which can be arbitrarily small, plus the optimization error $\delta$.

\subsection{Buffered Failure Probability Constraint}

We return to Example \ref{eBuffered} and the problem \eqref{eqn:bufferProbForm0} with a buffered failure probability constraint. As discussed in Subsection \ref{subsec:approx}, we aim to solve \eqref{eqn:bufferProbForm0} via the reformulation \eqref{eqn:bufferProbFormLagr} and the approximation \eqref{eqn:bufferProbFormLagrApprox}. Since the setting is more specific and also slightly different, we refine and adjust Theorem \ref{thm:optgap}.

The problem \eqref{eqn:bufferProbFormLagr} is defined using a space $\hat Z$, which is assumed to be a separable Banach space, but is augmented with $\reals^2$ to produce $Z = \hat Z \times \reals^2$. Using the product norm, this defines $(Z,\|\cdot\|_Z)$. While the underlying PDE only depends on $\hat z\in \hat Z$, we retain the notation $s:\Xi\times Z\to U$, which then involves a trivial extension. Likewise, the functions $\hat g_1$ and $\hat g_2$ from Example \ref{eBuffered} are extended to $U\times Z$. Thus, we are in the setting of Section \ref{sec:performance-bounds}, with $m=2$, $f$ given by $\hat f$ and $h(w) = w_1 + \iota_{\{0\}}(w_2)$.
The approximating problems have $h^\nu(w) = w_1 + y^\nu w_2 + \theta^\nu w_2^2$, for parameters $y^\nu\in\reals$ and $\theta^\nu \in (0,\infty)$, and
\begin{align*}
&g_1^\nu(u,z) = g_1(u,z)\\
&g_2^\nu(u,z) = \sigma + \gamma + \frac{1}{1-\alpha} \smax\Big(\hat g_2(u,z) - \gamma; ~\beta^\nu\Big), ~\mbox{ with } z = (\hat z, \gamma, \sigma).
\end{align*}

The constraint qualification in Assumption \ref{ass:CQ} can now be simplified.

\begin{assumption}{\rm (qualification for buffered constraint).}\label{ass:CQbuffer}
With the notation
\[
\psi(z) = \Ex\Big[\sigma + \gamma + \frac{1}{1-\alpha}\max\Big\{0, \hat g_2\big(s(\bfxi,z), z\big) - \gamma\Big\}\Big],~~~z = (\hat z, \gamma, \sigma),
\]
suppose that the following hold:

For fixed $n$ and  all $z_n\in A_n$, with $\psi(T_n(z_n))=0$, there exists $z^\nu_n\in A_n\to z_n$ as $\nu\to \infty$ such that $\psi(T_n(z_n^\nu))<0$ for all $\nu\in\nats$.

For all $z\in A$, with $\psi(z) = 0$, there exists $z^n\in A \cap Z^n\to z$ as $n\to \infty$ such that $\psi(z^n) < 0$ for $n\in\nats$.
\end{assumption}

\begin{proposition}{\rm (optimality gap under buffered probability constraint).}\label{cor:optgap2}
For fixed $\bar n\in\nats$, suppose that \eqref{eqn:bufferProbFormLagr} has a finite minimum value $\tau$, $A$ is closed, $f$ is continuous relative to $A$, the sample $\{\bfxi_1, \bfxi_2, \dots\}$ is iid as $\bfxi$, and, for all $n\geq \bar n$, $\inf \phi_n < \infty$, and suppose that Assumptions \ref{ass:rlsc}, \ref{ass:Gbound}, \ref{ass:upper}, \ref{ass:GboundInf}, and \ref{ass:CQbuffer} are satisfied. Moreover, let $\{y^\nu, \nu\in\nats\}$ be bounded, $\theta^\nu\to \infty$, $\beta^\nu\to 0$, and $\phi_n^\nu$ be the objective function in \eqref{eqn:bufferProbFormLagrApprox}.

Then, with probability one, the following hold: for any $\epsilon \in (0,\infty)$, there exists $n_\epsilon\geq \bar n$ such that if $n\geq n_\epsilon$, $\delta^\nu \to \delta \in [0,\infty)$, and $\{\bar z_n^\nu\in \delta^\nu\mbox{-}\nargmin \phi_n^\nu, \nu\in\nats\}$ has a cluster point $\bar z_n$, then $\bar z = T_n(\bar z_n) \in A$ and
\[
f(\bar z) + \Ex\Big[ \hat g_1\big(s(\bfxi, \bar z ), \bar z \big) \Big] \leq \tau + \epsilon + \delta, ~~~~~~\bprob\Big( \hat g_2\big(s(\bfxi, \bar z ), \bar z \big) > 0 \Big) \leq 1-\alpha.
\]
That is, $\bar z$ is feasible in \eqref{eqn:bufferProbForm0} with optimality gap of $\epsilon+\delta$.
\end{proposition}
\state Proof. We note that Assumption \ref{ass:smooth} holds automatically with $\epsilon^\nu = 2\beta^\nu/(1-\alpha)$; see the discussion in Subsection \ref{subsec:approx}.

Fix $n\geq \bar n$. Following the arguments in the proof of Theorems \ref{thm:perfguar} and \ref{thm:optgap}, we obtain that
\[
\nliminf_{\nu} \phi_n^\nu(z_n^\nu) \geq \phi_n(z_n) \mbox{ whenever } z_n^\nu \to z_n \mbox{ as } \nu\to \infty
\]
provided that $\nliminf_\nu h^\nu(w^\nu) \geq h(w)>-\infty$ for $w^\nu\to w=(w_1,w_2)$. If $w_2 \neq 0$, then $h(w) = \infty$ and $h^\nu(w^\nu) \to \infty$ because $\{y^\nu, \nu\in\nats\}$ is bounded and $\theta^\nu\to \infty$. If $w_2 = 0$, then $h(w) = w_1$ and $h^\nu(w^\nu) \geq w_1^\nu + y^\nu w_2^\nu \to w_1$  after again using the boundedness of $y^\nu$. Thus, the liminf-requirement holds.

We next show that for each $\breve z_n\in A_n$, there exists $\breve z_n^\nu\in A_n\to \breve z_n$ such that
	\begin{equation}\label{eqn:bufferproof0}
	\nlimsup_{\nu} \phi_n^\nu(\breve z_n^\nu) \leq \phi_n(\breve z_n).
	\end{equation}
This holds trivially if $\breve z_n \not\in A_n$. For $\breve z_n = (\hat z_n, \gamma, \sigma)\in A_n$, we argue as follows. Let $\{\xi_1, \xi_2, \dots\}$ be an event for which the convergence in \eqref{eqn:convprop} holds.
As in the proof of Theorem \ref{thm:perfguar}, we know that
\begin{equation}\label{eqn:convsampl3}
\frac{1}{\nu} \sum_{j=1}^{\nu} G_n^\nu(\xi_j,z_n^\nu) \to \Ex\big[G_n(\bfxi,z_n)\big]~ \mbox{ whenever } ~ z_n^\nu \in A_n \to z_n.
\end{equation}
We consider two cases. Let $\psi$ be as defined in Assumption \ref{ass:CQbuffer}.

First, suppose that $\psi(T_n(\breve z_n)) \neq 0$. Then, we set $\breve z_n^\nu = \breve z_n$ for all $\nu$ and obtain \eqref{eqn:bufferproof0} because the right-hand side equals infinity.

Second, suppose that $\psi(T_n(\breve z_n)) = 0$. By Assumption \ref{ass:CQbuffer}, there exist $z_n^k = (\bar z_n^k,\gamma^k,\sigma^k) \in A_n\to \breve z_n$ as $k\to \infty$ such that $\psi(T_n(z_n^k)) < 0$.
Thus, for each $k$, there exists by \eqref{eqn:convsampl3} $\nu_k$ such that for all $\nu\geq \nu_k$ one has
\[
\sigma^k + \gamma^k + \frac{1}{1-\alpha}\frac{1}{\nu} \sum_{j=1}^\nu \smax\bigg(\hat g_2\Big(s\big(\xi_j,T_n(z_n^k) \big),  T_n(z_n^k) \Big)  - \gamma^k; ~\beta^\nu \bigg) \leq \frac{1}{2}\psi\big(T_n(z_n^k) \big).
\]
This inequality means that, for each $\nu$, there is $k(\nu)$ such that
\[
\sigma^{k(\nu)} + \gamma^{k(\nu)} + \frac{1}{1-\alpha}\frac{1}{\nu} \sum_{j=1}^\nu \smax\bigg(\hat g_2\Big(s\big(\xi_j,T_n(z_n^{k(\nu)}) \big),  T_n(z_n^{k(\nu)}) \Big)  - \gamma^{k(\nu)}; ~\beta^\nu \bigg) \leq 0
\]
and this can be done such that $k(\nu)\to \infty$ as $\nu\to \infty$. Set $\tilde \sigma^\nu \geq 0$ equal to the negative of the left-hand side in this inequality. Thus, $\breve z_n^\nu = (\bar z_n^{k(\nu)}, \gamma^{k(\nu)}, \sigma^{k(\nu)}+\tilde\sigma^\nu)\in A_n$ and tends to $\breve z_n$ as $\nu\to \infty$.
Moreover,
\begin{align*}
h^\nu\Bigg( \frac{1}{\nu} \sum_{j=1}^{\nu} G_n^\nu(\xi_j,\breve z_n^\nu) \Bigg) & = \frac{1}{\nu} \sum_{j=1}^{\nu} \hat g_1\Big(s\big(\xi_j, T_n(\breve z_n^\nu) \big), T_n(\breve z_n^\nu) \Big) + y^\nu\cdot 0 + \theta^\nu \cdot 0^2\\
& ~~~\to \Ex\Big[ \hat g_1\big(s(\bfxi, T_n(\breve z_n) ), T_n(\breve z_n) \big) \Big] = h\Big(\Ex\big[G_n(\bfxi,\breve z_n)\big]\Big).
\end{align*}
Since $\breve z_n^\nu\in A_n$ and $f_n$ is continuous on $A_n$, one obtains $\phi_n^\nu(\breve z_n^\nu) \to \phi_n(\breve z_n)$. Consequently, with probability one, we have established that $\phi_n^\nu \eto \phi_n$ as $\nu\to \infty$; see Definition \ref{def:epiconv}.

Next, we confirm that $\phi + \iota_{Z^n} \eto \phi$ as $n\to \infty$. Let $z^n\to z$. Since $h$ is lsc, we can argue as in the proof of Theorem \ref{thm:optgap} and conclude that
\[
\nliminf_n \big(\phi(z^n) + \iota_{Z^n}(z^n)\big) \geq \phi(z).
\]
Let $z = (\hat z, \gamma, \sigma)\in Z$. We also need to construct $z^n\to z$ such that
\begin{equation}\label{eqn:limsupcondMain2}
	\nlimsup_n \big(\phi(z^n) + \iota_{Z^n}(z^n)\big) \leq \phi(z).
\end{equation}
We can assume that $z\in A$ because otherwise the right-hand side will be infinity. If $\psi(z) \neq 0$, then $\phi(z) = \infty$ and \eqref{eqn:limsupcondMain2} holds trivially again. If $\psi(z) = 0$, then, by the second part of Assumption \ref{ass:CQbuffer}, there is $(\hat z^n, \gamma^n, \sigma^n) \in  A \cap Z^n\to z$ such that $\psi((\hat z^n, \gamma^n, \sigma^n)) < 0$. Set $\tilde \sigma^n = -\psi((\hat z^n, \gamma^n, \sigma^n))$, which then vanishes as $n\to \infty$, and also construct $z^n = (\hat z^n, \gamma^n, \sigma^n+\tilde \sigma^n)$. Since $\psi(z^n) = 0$, we conclude that
\[
\phi(z^n) + \iota_{Z^n}(z^n) = f(z^n) + \Ex\Big[ \hat g_1\big(s(\bfxi, z^n ), z^n \big) \Big] \to f(z) + \Ex\Big[ \hat g_1\big(s(\bfxi, z ), z \big) \Big] = \phi(z)
\]
after invoking Lemma \ref{lem:llnInfDim}.
We have established that $\phi + \iota_{Z^n} \eto \phi$; cf. Definition \ref{def:epiconv}. The assumptions of Theorem \ref{thm:minval} therefore holds and $\phi( T_n(\bar z_n) ) \leq \inf \phi + \epsilon + \delta$. Since $\inf \phi = \tau$ is finite, $T_n(\bar z_n)$ is feasible in \eqref{eqn:bufferProbFormLagr}. The conclusion then follows by reversing the arguments leading from \eqref{eqn:bufferProbForm0} to \eqref{eqn:bufferProbForm} and to \eqref{eqn:bufferProbFormLagr}.\eop

The proposition shows that a cluster point constructed by the approximating problems \eqref{eqn:bufferProbFormLagrApprox} produces a feasible control for the actual problem that is suboptimal with tolerance $\epsilon + \delta$. There is much flexibility in choosing the multipliers $y^\nu$. A penalty method remains possible by setting  $y^\nu = 0$, but update rules for $y^\nu$ from augmented Lagrangian methods may be computationally superior. Under additional assumptions, we conjecture that ``exactness'' is possible and then $\theta^\nu$ can remain bounded; we omit the details and refer to \cite[Section 6.B]{primer}.

\section{Application}\label{sec:applications}

This section examines the assumptions of Section \ref{sec:performance-bounds} in the context of a distributed optimal control problems. Specifically, for $d =1, 2$, or $3$, let $D \subset \bR^d$ be a bounded open set representing a physical domain and let $\partial D$ be the corresponding boundary. Set $\bar D = D \cup \partial D$. We denote by $L^2(D)$ and $L^2(\partial D)$ the spaces of Lebesgue square-integrable functions from $D$ and $\partial D$ to $\reals$, respectively. As before, $(\Omega, \cF, \bP)$ is a probability space. We consider thermal conduction modeled by the following PDE with random conductivity coefficient: For a.e.\ $\omega \in \Omega$,
	\begin{equation}\label{eq:heat}
		\begin{split}
			-\nabla \cdot (\bfxi(\omega) \nabla u) & = c_1 z \quad \text{in } D\\
			\bfxi(\omega) \nabla u \cdot \vec{n} & = c_2 (s_e - u ) \quad \text{on } \partial D,
		\end{split}
	\end{equation}
where the state variable $u: \bar{D} \to \bR$ represents the unknown temperature field, $s_e\in L^2(\partial D)$ is a given exterior temperature along the boundary, $z \in Z = L^2(D)$ is the control in a nonempty admissible set $A = \{z \in Z~|~\underline{z} \leq z(x) \leq \overline{z} \mbox{ for a.e. } x\in D\}$, $c_1: D \to [0, \infty)$ and $c_2:\partial D\to [0,\infty)$ are known coefficients that are uniformly bounded from above, $\vec{n}$ is a unit normal vector pointing outside $\bar D$, and $\bfxi(\omega):\bar D\to (0,\infty)$ is a thermal conductivity field. Specifically,
\[
	\bsxi(\omega)(x) = \exp \bigg(b_0(x) + \sum_{j = 1}^J b_j(x) \bfy_j(\omega)\bigg),	
\]
where $J \in \bN$,  $\bfy_1, \dots, \bfy_J$ are iid random variables defined on the probability space $(\Omega, \cF, \bP)$ and $b_0, \dots, b_J:\bar D\to \reals$ are given functions. This induces a probability space $(\Xi,\cB,P)$ as discussed in Section \ref{sec:formulation}, where elements $\xi\in \Xi$ are functions from $\bar D$ to $(0,\infty)$. One example of such a random field has $\log(\bsxi)$ given by a truncated Kahunen--Lo\`eve expansion of a Gaussian random field $\cN(b_0, \cC)$ with mean $b_0$ and covariance $\cC$ \cite{SchwabTodor06}, where $b_j = \sqrt{\lambda_j} \phi_j$, $(\lambda_j, \phi_j)$ is an eigenpair of $\cC$ and $\bfy_1, \dots, \bfy_J$ follow the standard normal distribution. In another example, $\bsxi$ represents a piecewise random thermal conductivity coefficient with $b_j$ denoting a rescaled characteristic function supported on a subdomain $D_j \subset \bar D$.

\begin{assumption}\label{ass:xi-ass}
There exist $\underline{c},\overline{c}: \Xi \to  (0, \infty)$ such that
	\begin{equation}\label{eq:xi-bound}
		0  < \underline{c}(\xi) = \essinf_{x\in \bar D} \xi(x) \leq \esssup_{x\in \bar D}\xi (x) = \overline{c}(\xi) < \infty, \; \text{for a.e. } \xi \in \Xi
	\end{equation}
and, regardless of $0 \leq p \leq 1$ and $0 \leq q \leq 3$, one has
                \begin{equation}\label{eq:integrable-xi}
                                \int_\Xi \frac{(\overline{c}(\xi))^p}{(\underline{c}(\xi))^q} dP(\xi) < \infty.
                \end{equation}
\end{assumption}

Following the notation and terminology from \eqref{eqn:actualpde0}, we consider $U = L^2(D)$ and two performance functions: the discrepancy between the solution and a desired temperature field $s_d: D \to \reals$ and the shortfall of the solution average over $D_t\subset D$ relative to a threshold temperature $s_t$. They are captured by $g_1,g_2:U \times Z \to \reals$ with
\begin{equation}\label{eq:g1g2}
	g_1(u,z) = \int_D \big(u(x) - s_d(x)\big)^2 dx, ~~~~~ g_2(u,z) = s_t - \int_{D_t} u(x) dx.
\end{equation}
We adopt the cost function $f(z) = \theta \|z\|_{Z}^2$ for some $\theta > 0$, which is continuous, and thus satisfying the continuity assumption of $f$ in Theorem \ref{thm:optgap} and Proposition \ref{cor:optgap2}.

Let $V = H^1(D) = \{v\in L^2(D)~|~ |\nabla v| \in L^2(D)\}$ with the norm $\|v\|_V = \|v\|_{L^2(D)} + \||\nabla v|\|_{L^2(D)}$. For fixed $\xi \in \Xi$ and $z\in Z$, $u\in V$ is a weak solution of \eqref{eq:heat} if
\begin{equation}\label{eq:heat-weak}
	a(u, v; \xi) = \ell(v) + b(z, v) \quad \forall v \in V,
\end{equation}
where
\begin{align*}
a(u, v;\xi) = \int_D \xi(x) \nabla u(x) \cdot \nabla v(x) dx + \int_{\partial D} c_2(x) u(x) v(x) dx \quad \forall u, v \in V, \\
~ b(z, v) = \int_D c_1(x) z(x) v(x) dx, ~ \text{ and } \ell(v) = \int_{\partial D} c_2(x) s_e(x) v(x) dx \quad \forall v \in V.
\end{align*}
A solution $u$ of \eqref{eq:heat-weak} defines a mapping from $\Xi\times Z$ to $V$. Since the performance functions $g_1,g_2$ are defined on $L^2(D) \times L^2(D)$, it becomes more natural to adopt $U = L^2(D)$ as the range space for the solution mapping. Thus, the solution mapping $s:\Xi\times Z\to U$ is given as $s(\xi,z) = u$, where $u$ is the solution of \eqref{eq:heat-weak} under input $\xi \in \Xi$ and $z \in Z$. We are in the setting of the earlier sections with $Z = U = L^2(D)$.

\subsection{Approximations}

We use a finite element method to approximate solutions and controls. Let $X_h^p$ denote a finite element space of piecewise polynomials of total degree $p\in \nats\cup \{0\}$ on each element $K$ of a triangulation $\cT_h$ with mesh size $h$. For simplicity, we set $Z^n = X_{h}^{0}$ for some $h$ such that $n$ represents the number of elements. Let $T_n:\reals^n \to Z^n$ be given by
\begin{equation}\label{eq:z-approx}
z^n =  T_n(z_n) = \sum_{i = 1}^n z_n^i \psi_i(\cdot)~~~~~~\mbox{ for } z_n = (z_n^1, z_n^2, \dots, z_n^n)\in\reals^n,
\end{equation}
where $\psi_i(x) = 1$ for $x$ in element $K_i$ and zero otherwise. Let $V^\nu = X_h^1$ for some $h$ depending on an index $\nu$, with $r^\nu$ representing the number of finite element basis functions in $V^\nu$, i.e., the number of degrees of freedom. Although not common in practice, $Z^n$ and $V^\nu$ can be constructed with different meshes. The Galerkin approximation of the problem \eqref{eq:heat-weak} reads: for $\xi \in \Xi$ and control $z^n \in Z^n$, find $u^\nu \in V^\nu$ such that
\begin{equation}\label{eq:heat-weak-h}
		a(u^\nu, v^\nu; \xi) = \ell(v^\nu) + b(z^n, v^\nu) \quad \forall v^\nu \in V^\nu.
\end{equation}
We remark that \eqref{eq:heat-weak-h} can be extended to any control $z\in Z$.
The approximating solution mapping $s^\nu:\Xi\times Z \to U$ is therefore defined as $s^\nu(\xi,z) = u^\nu$, where $u^\nu$ is the solution of \eqref{eq:heat-weak-h} under input $\xi \in \Xi$ and $z \in Z$.

These approximations have the following properties.

\begin{proposition}\label{prop:Tnzn}
	The mapping $T_n: \bR^n \to Z^n$ in \eqref{eq:z-approx} is continuous, with
	\begin{equation}\label{eq:z-n-bound}
		\big\|T_n(z_n)\big\|_Z  \leq {c}_K^{1/2} h^{d/2} \|z_n\|_2 \quad \forall z_n\in\reals^n,
	\end{equation}
	for a constant ${c}_K < \infty$ independent of $h$, $n$, and $z_n$. Moreover, for any $z \in A$, there exists $\{z_n\in\reals^n, n\in\nats\}$ such that
	\begin{equation}
		\big\|z - T_n(z_n)\big\|_Z \to 0 ~\text{ as } ~n \to \infty.
	\end{equation}
\end{proposition}
\state Proof. By definition of the piecewise constant finite element approximation, we have
\begin{equation}
	\begin{split}
	\big\|T_n(z_n)\big\|_{Z}^2 & = \int_{D} \big(T_n(z_n)(x)\big)^2 dx = \sum_{i= 1}^n \sum_{j = 1}^n z_n^i z_n^j \int_{D} \psi_i(x) \psi_j(x) dx \\
	& = \sum_{i = 1}^n (z_n^i)^2 |K_i| \leq c_K h^{d} \|z_n\|_2^2,
	\end{split}
\end{equation}
where the inequality holds because each element $K_i \in \cT_h$ of mesh size $h$, has length, area, or volume  $|K_i| \leq c_K h^d$ in $d = 1, 2, 3$, respectively, for some constant $c_K < \infty$ independent of $n$. Moreover, for any function $z \in A$, which is measurable, bounded, and supported in the bounded domain $D$, the simple function approximation $T_n(z_n)$ is uniformly convergent. \eop

\begin{proposition}\label{prop:solutions}
Suppose that Assumption \ref{ass:xi-ass} holds. For $\xi \in \Xi$ and $z\in Z$, there exist a unique solution $s(\xi, z) \in V$ of \eqref{eq:heat-weak} and a unique solution $s^\nu(\xi, z) \in V^\nu \subset V$ of \eqref{eq:heat-weak-h}. Moreover, we have the stability estimate
		\begin{equation}\label{eq:stability}
			\big\|s(\xi, z)\big\|_{V}, \big\|s^\nu(\xi, z)\big\|_{V}  \leq c^s_2(\xi) + c^s_1(\xi) \|z\|_{Z},
		\end{equation}
		where $c^s_2$ and $c^s_1$ are given in \eqref{eq:lambdas}. Furthermore, there holds the continuity estimate
		\begin{equation}\label{eq:s-continuity}
			\big\|s(\xi, z) - s(\xi, z')\big\|_{V}, \big\|s^\nu(\xi, z) - s^\nu(\xi, z')\big\|_{V} \leq  c_z^s(\xi) \|z - z'\|_{Z},
		\end{equation}
		where $c_z^s(\xi)$ is given in \eqref{eq:c-c-s}, and the continuity estimate
		\begin{equation}\label{eq:s-continuity-xi}
			\big\|s(\xi, z) - s(\xi', z)\big\|_V, \big\|s^\nu(\xi, z) - s^\nu(\xi', z)\big\|_V \leq  c_\xi^s(\xi,z) \|\xi - \xi'\|_{L^\infty(D)},
		\end{equation}
		where $c_\xi^s(\xi, z)$ is given in \eqref{eq:c-xi-s}. Since $V = H^1(D) \subset U = L^2(D)$ and $\|u\|_U \leq \|u\|_V$ for any $u \in V$, all the estimates above also hold in $\|\cdot\|_U$.

		Finally, we have the error estimate
		\begin{equation}\label{eq:error-estimate}
			\big\|s(\xi, z) - s^\nu(\xi, z)\big\|_{U} \leq c_\nu^s(\xi)\big\|s(\xi, z)\big\|_{V} (r^\nu)^{-1/d} \to 0 \text{ as } \nu \to \infty,
		\end{equation}
		where $c_\nu^s(\xi)$ is given in \eqref{eq:c-nu-s}, $r^\nu$ is the number of degrees of freedom of $V^\nu$, which satisfies $r^\nu \to \infty$ as $\nu \to \infty$, and $d$ is the dimension of the physical domain $D \subset \bR^d$.
	\end{proposition}
\state Proof. The proof follows that for finite element Galerkin approximation of deterministic elliptic PDE (see, for example, \cite[\S 4]{QuarteroniQuarteroni09}), but is enriched with additional details to facilitate verification of Assumptions \ref{ass:Gbound} and \ref{ass:upper}.

By \eqref{eq:xi-bound} and the Poincar\'e inequality $\|u\|_{L^2(D)}  \leq c_P \|\nabla u\|_{L^2(D)}$, with the Poincar\'e constant $c_P > 0$, one obtains
	\begin{equation}\label{eq:coercivity}
		a(u, u;\xi) \geq \underline{c}(\xi) \|\nabla u\|^2_{L^2(D)} \geq \frac{\underline{c}(\xi)}{\sqrt{1+c_P}} \|u\|_V^2, \quad \forall u \in V.
	\end{equation}
The Cauchy--Schwarz inequality, \eqref{eq:xi-bound}, $\overline{c}_2 = \|c_2\|_{L^\infty(D)}$,  and the trace inequality $\|v\|_{L^2(\partial D)} \leq c_T \|v\|_V$ with constant $c_T > 0$ produce
	\begin{equation}
		a(u, v; \xi) \leq \big(\overline{c}(\xi) + \overline{c}_2 c_T^2\big)\|u\|_V \|v\|_V.
	\end{equation}
Similarly, with $\overline{c}_1 = \|c_1\|_{L^\infty(D)}$ and $\|v\|_{L^2(D)} \leq \|v\|_V$, we have
	\begin{equation}\label{eq:b-bound}
		b(z, v) \leq \overline{c}_1 \|z\|_Z \|v\|_V, \quad \forall z \in Z, v \in V.
	\end{equation}
Again, by the Cauchy--Schwarz inequality and the trace inequality, we reach
	\begin{equation}
		\ell(v) \leq \overline{c}_2 c_T \|s_e\|_{L^2(\partial D)} \|v\|_V, \quad \forall v \in V.
	\end{equation}
Combining the above four inequalities, we conclude that there exists a unique solution $s(\xi, z) \in V$ of \eqref{eq:heat-weak} by the Lax--Milgram theorem and also a unique solution $s^\nu(\xi, z) \in V^\nu$ of \eqref{eq:heat-weak-h} since $V^\nu \subset V$, for which the above four inequalities also hold. Moreover, replacing $v = u$ in \eqref{eq:heat-weak} and $v^\nu = u^\nu$ in \eqref{eq:heat-weak-h}, by the above inequalities, we have the stability estimate \eqref{eq:stability},
	with $c^s_2$ and $c^s_1$ given by
	\begin{equation}\label{eq:lambdas}
		c^s_2(\xi) = \frac{\overline{c}_2 c_T \|s_e\|_{L^2(\partial D)} \sqrt{1+c_P}}{\underline{c}(\xi)} \quad \text{ and } \quad c^s_1(\xi) = \frac{\overline{c}_1 \sqrt{1+c_P} }{\underline{c}(\xi)}.
	\end{equation}

	For the continuity of $s(\xi, \cdot)$, we take $z, z' \in Z$ in \eqref{eq:heat-weak} and denote $s = s(\xi, z)$ and $s' = s(\xi, z')$ as the corresponding solutions, respectively. By subtracting the two equations \eqref{eq:heat-weak} at $z, z' \in Z$, we have
	\begin{equation}
		a(s - s', v; \xi) = b(z-z', v) \quad \forall v \in V.
	\end{equation}
	By $v = s-s'$, we obtain the continuity estimate \eqref{eq:s-continuity} using \eqref{eq:coercivity} and \eqref{eq:b-bound}, with % $c_z^s(\xi)$ given by
	\begin{equation}\label{eq:c-c-s}
		c_z^s(\xi) = \frac{\overline{c}_1 \sqrt{1+c_P}}{\underline{c}(\xi)}.
	\end{equation}
	By the same argument, the same continuity estimate holds for $s^\nu(\xi, \cdot)$.

	Similarly, for the continuity of $s(\cdot, z)$, we take $\xi, \xi' \in \Xi$ in \eqref{eq:heat-weak} and denote $s = s(\xi, z)$ and $s' = s(\xi', z)$ as the corresponding solutions, respectively. By subtracting the two equations \eqref{eq:heat-weak} at $\xi, \xi' \in \Xi$, we have
	\begin{equation}
		\int_D \xi(x) \nabla \big(s(x) - s'(x)\big) \cdot \nabla v(x) dx = \int_D \big(\xi'(x) - \xi(x)\big) \nabla s'(x) \cdot \nabla v(x) dx.
	\end{equation}
	By taking $v = s - s'$, we obtain the continuity estimate \eqref{eq:s-continuity-xi} using \eqref{eq:coercivity} by noting that $c_2 \geq 0$ and the stability estimate \eqref{eq:stability}, with $c_\xi^s(\xi, z)$ given by
	\begin{equation}\label{eq:c-xi-s}
		c_\xi^s(\xi,z) = \frac{\sqrt{1+c_P}}{\underline{c}(\xi)} \big( c^s_2(\xi) + c^s_1(\xi) \|z\|_Z\big).
	\end{equation}
	The same holds for the approximation solution $s^\nu$ by the same argument.

	The error estimate \eqref{eq:error-estimate} follows the proof of \cite[\S 4, Thm. 4.7]{QuarteroniQuarteroni09}, which satisfies
	\begin{equation}
		\big\|s(\xi, z) - s^\nu(\xi, z)\big\|_{U} \leq h \hat c \frac{\overline{c}(\omega)}{\underline{c}(\xi)}\big\|s(\xi, z)\big\|_{V},
	\end{equation}		
		where $h$ is the finite element mesh size for the discretization of the domain $D$ and $\hat c \in (0, \infty)$ is a constant, independent of $\xi$ and $z$.
	Note that the mesh size $h$ is related to the number of degrees of freedom $\nu$ as $h = O((r^\nu)^{-1/d})$ in $D \subset \bR^d$, i.e., $h \leq \breve{c} (r^\nu)^{-1/d}$	for some $\breve{c}\in (0,  \infty)$. The mesh size is thus tending to zero. Then, we obtain \eqref{eq:error-estimate} with $c_\nu^s(\xi)$ given by
	\begin{equation}\label{eq:c-nu-s}
		c_\nu^s(\xi) = \hat c \breve{c} \overline{c}(\xi)/\underline{c}(\xi),
	\end{equation}
	which concludes all the estimates. \eop

To address Assumption \ref{ass:rlsc}, we note that $s(\xi, \cdot): Z \to U$ is continuous due to \eqref{eq:s-continuity} and $s(\cdot,z)$ and $s^\nu(\cdot,z)$ are measurable due to \eqref{eq:s-continuity-xi}.

	\begin{proposition}\label{prop:g1g2}
		The performance functions $g_1$ and $g_2$ defined in \eqref{eq:g1g2} satisfy the continuity estimates: for all $u, u' \in U, z, z' \in Z$,
		\begin{equation}
			\begin{split}
				\big|g_1(u, z) - g_1(u', z')\big| & \leq \big(\|u\|_{U} + \|u'\|_{U} + 2\|s_d\|_{U}\big) \|u - u'\|_{U}\\
				\big|g_2(u, z) - g_2(u', z')\big| & \leq |D_t| \|u - u'\|_{U},
			\end{split}
		\end{equation}
		where $|D_t|$ is the size of the subdomain $D_t \subset D$. This implies  Assumption \ref{ass:rlsc} and Assumption \ref{ass:upper} for $g_1$ and $g_2$. Moreover, for all $u \in U, z \in Z$, the performance functions are bounded from above by
\[
\big|g_1(u, z)\big| \leq \big(\|u\|_{U} + \|s_d\|_{U}\big)^2,~~~~~~~				\big|g_2(u, z)\big| \leq |s_t| + |D_t| \|u\|_{U}.
\]
\end{proposition}
\state Proof. By definition of $g_1$ in \eqref{eq:g1g2}, we have
	\begin{equation*}
		\begin{split}
		&	\big|g_1(u, z) - g_1(u', z')\big| = \left| \int_D \big(u(x) - s_d(x)\big)^2 - \big(u'(x) - s_d(x)\big)^2 dx\right|\\
			&  \leq \big(\|u\|_{U} + \|u'\|_{U} + 2\|s_d\|_{U}\big) \|u  - u'\|_{U},
		\end{split}
	\end{equation*}
	which concludes the continuity estimate. The upper bound for $g_1$ is satisfied by the Cauchy-Schwarz inequality. For $g_2$, we have
	\begin{equation*}
		\big|g_2(u, z) - g_2(u', z')\big| = \left| \int_{D_t} \big(u(x) - u'(x)\big) dx \right| \leq |D_t| \; \|u - u'\|_{U},		
	\end{equation*}
	which concludes the continuity estimate. The upper bound for $g_2$ is also satisfied by the Cauchy-Schwarz inequality. \eop

	Note that the continuity of $g_1$ and $g_2$ implies the measurability of $g_1(s(\xi, z), z)$ and $g_2(s(\xi, z), z)$ w.r.t.\ $\xi$, by the measurability of the solution $s(\xi, z)$ w.r.t.\ $\xi$ from Proposition \ref{prop:solutions} and the Doob--Dynkin lemma \cite[Lemma 2.12]{Oksendal03}, which verifies Assumption \ref{ass:rlsc} on the measurability of $G(\xi, z) = (g_1(s(\xi, z), z), g_2(s(\xi, z), z))$ w.r.t.\ $\xi$.

	\begin{proposition}
		Let $n\in \nats$ be fixed and consider the setting of this section under Assumption \ref{ass:xi-ass}. For $\bar{z}_n \in A_n$, $\rho > 0$, and $z_n \in \ball(\bar{z}_n, \rho) \cap A_n$, let $z^n = T_n(z_n) \in Z^n$. For $i = 1, 2$, there exists an integrable function $\kappa_i: \Xi \to \bR$, such that
		\begin{equation*}
				\Big|g_i\big(s(\xi, z^n), z^n\big)\Big| \leq \kappa_i(\xi),
		\end{equation*}
		which also hold when $z^n$ is replaced by any $z \in A$ such that $\|z - \bar{z}\|_Z \leq \rho$ for $\bar{z} \in A$.
		Moreover, there exist integrable $\lambda_i: \Xi \to \bR$ and $\Delta_i^\nu(\bar{z}_n) \to 0$ as $\nu \to \infty$, such that
		\begin{equation*}
			\Big|g_i\big(s(\xi, z^n), z^n\big) - g_i\big(s^\nu(\xi, z^n), z^n\big)\Big| \leq \lambda_i(\xi) \Delta_i^\nu(\bar{z}_n),
		\end{equation*}
		where $\kappa_i, \lambda_i, \Delta_i^\nu$ are defined in the proof.
	\end{proposition}
	\state Proof. For $z_n \in \ball(\bar{z}_n, \rho) \cap A_n$, one has $\|z^n\|_Z = \|T_n(z_n)\|_Z < \infty$ by Proposition \ref{prop:Tnzn}.

	First, we prove the upper bound and the error estimate for $g_1$. By definition and Proposition \ref{prop:g1g2}, we have
\[
g_1\big(s(\xi, z^n), z^n\big)  \leq \big(\|s(\xi, z^n)\|_U + \|s_d\|_{U}\big)^2 \leq  \big(c^s_2(\xi) + c^s_1(\xi) \|z^n\|_Z + \|s_d\|_{U}\big)^2,
\]
where we used the stability estimate \eqref{eq:stability}. The right-hand side expression furnishes a value for $\kappa_1(\xi)$. The random variable $\kappa_1$ is integrable if $(c^s_2(\cdot))^2$, $(c^s_1(\cdot))^2$, $c^s_2(\cdot)c^s_1(\cdot)$, $c^s_2(\cdot)$, $c^s_2(\cdot)$ are all integrable. By definition of $c^s_2(\xi)$ and $c^s_2(\xi)$ in \eqref{eq:lambdas}, this is satisfied because $1/\underline{c}(\cdot)$ and $1/(\underline{c}(\cdot))^2$ are integrable by Assumption \ref{ass:xi-ass}. Moreover, the integrability of $\kappa_1$ holds when $z^n \in Z^n$ is replaced by any $z \in Z$ since $\|z\|_Z < \infty$.
	
By the continuity estimate for $g_1$ in Proposition \ref{prop:g1g2}, and the error and stability estimates in Proposition \ref{prop:solutions} for $s(\xi,z^n)$ and $s^\nu(\xi,z^n)$, we have
		\begin{equation*}
		\begin{split}
			& \Big|g_1\big(s(\xi,z^n), z^n\big) - g_1\big(s^\nu(\xi,z^n), z^n\big)\Big|\\
			& \leq c_\nu^s(\xi) \left(\big\|s(\xi,z^n)\big\|_U + \big\|s^\nu(\xi,z^n) \big\|_U + 2\|s_d\|_{U}\right)  \|s\|_{V} (r^\nu)^{-1/d} \\
			& \leq c_\nu^s(\xi) \Big(
			(c_2^s(\xi))^2	+ 2c_2^s(\xi)c_1^s(\xi) \|z^n\|_Z + (c_1^s(\xi)\|z^n\|_Z)^2	\\
			& ~~~ + c_2^s(\xi) \|s_d\|_{U} + c_1^s(\xi) \|s_d\|_{U} \|z^n\|_Z \Big)  (r^\nu)^{-1/d} \leq \lambda_1(\xi) \Delta^\nu_1(\bar z_n),
		\end{split}
	\end{equation*}
	where
		\begin{align*}
			\lambda_1(\xi) & = c_\nu^s(\xi) \Big((c_2^s(\xi))^2	+ 2c_2^s(\xi)c_1^s(\xi)  + (c_1^s(\xi))^2	
			 + c_2^s(\xi) + c_1^s(\xi)\Big)\\
			\Delta^\nu_1(\bar z_n) & = \hspace{-0.06cm} \sup_{z_n \in \ball(\bar{z}_n, \rho)} \hspace{-0.2cm}\big( 1 + \|T_n(z_n)\|_Z + \|T_n(z_n)\|_Z^2  + \|s_d\|_{U} + \|s_d\|_{U} \|T_n(z_n)\|_Z \big) (r^\nu)^{-1/d}
	\end{align*}
	which satisfies $\Delta^\nu_1(\bar z_n) \to 0$ as $\nu \to \infty$.
	By the definition of $c_1^s$ and $c_2^s$ in \eqref{eq:lambdas}, we have that $\lambda_1$ is integrable, as long as $\overline{c}(\cdot)/(\underline{c}(\cdot))^3$ and $\overline{c}(\cdot)/(\underline{c}(\cdot))^2$ are integrable, which is the case by Assumption \ref{ass:xi-ass}.

Second, for $g_2$, by Proposition \ref{prop:g1g2} we have
\[
\Big|g_2\big(s(\xi, z^n), z^n\big)\Big| \leq |s_t| + |D_t| \big\|s(\xi, z^n)\big\|_U \leq |s_t| + |D_t| \big(c^s_2(\xi) + c^s_1(\xi) \|z^n\|_Z\big),
\]
where we used the stability estimate \eqref{eq:stability}, with $|D_t|$ measuring the size of $D_t$. The right-hand side furnishes an expression for $\kappa_2(\xi)$. This defines an integrable random variable $\kappa_2$ because $c^s_2(\cdot)$ and $c^s_2(\cdot)$ in \eqref{eq:lambdas} are integrable as a result of that $1/\underline{c}(\cdot)$ is integrable by Assumption \ref{ass:xi-ass}. Moreover, the integrability of $\kappa_2$ holds when $z^n \in Z^n$ is replaced by $z \in Z$ as $\|z\|_Z < \infty$.

	By the continuity estimate for $g_2$ in Proposition \ref{prop:g1g2}, and the error and stability estimates in Proposition \ref{prop:solutions} for $s(\xi,z^n)$ and $s^\nu(\xi,z^n)$, we have
	\begin{equation*}
		\begin{split}
				\Big|g_2\big(s(\xi,z^n), z^n\big) - g_2\big(s^\nu(\xi,z^n), z^n\big)\Big| & \leq |D_t| c_\nu^s(\xi) \big(c^s_2(\xi) + c^s_1(\xi) \|z^n\|_Z\big) (r^\nu)^{-1/d} \\
				& \leq \lambda_2(\xi)  \Delta_2^\nu(\bar{z}_n),
		\end{split}
	\end{equation*}
	where $\lambda_2(\xi) = |D_t| c_\nu^s(\xi) (c^s_2(\xi) + c^s_1(\xi))$,	which is integrable as long as $\overline{c}(\cdot)/(\underline{c}(\cdot))^2$ is integrable, which holds by Assumption \ref{ass:xi-ass}. Moreover,
	\begin{equation*}
		\Delta_2^\nu(\bar{z}_n) = \sup_{z_n \in \ball(\bar{z}_n, \rho)} \big(1 + \|T_n(z_n)\|_Z\big) (r^\nu)^{-1/d},
	\end{equation*}
	which satisfies $\Delta_2^\nu(\bar{z}_n)  \to 0 $ as $\nu \to \infty$. % Moreover, the measurability of $g_2$ in Assumption \ref{ass:rlsc} is satisfied by following the same argument for $g_1$.
	\eop

To this end, all the assumptions made in Section \ref{sec:intermediate} are satisfied under Assumption \ref{ass:xi-ass} for the example presented in this section.

	\bibliographystyle{siamplain}
	% \bibliography{refs,references}

\end{document}